\documentclass[11pt]{amsart}
\usepackage{amsmath, amsfonts, latexsym, amssymb, graphicx, mathtools,  wasysym, enumitem, pifont}
\usepackage{epsfig}
\usepackage{amssymb, amsmath}
\def\cat{{\text CAT}}
\newtheorem{thm}{Theorem}[section]
\newtheorem{defn}[thm]{Definition}
\newtheorem{proposition}[thm]{Proposition}
\newtheorem{lemma}[thm]{Lemma}

\begin{document}

\title[Local connectivity of right-angled Coxeter group boundaries]
{Local connectivity of right-angled Coxeter group boundaries}

\author[Michael Mihalik]{Michael Mihalik}
\address{Department of Mathematics\\
        Vanderbilt University\\
        Nashville, TN 32340}
\email{mihalik@math.vanderbilt.edu}

\author[Kim Ruane]{Kim Ruane}
\address{Department of Mathematics\\
        Tufts University\\
        Medford MA 02155}
\email{kim.ruane@tufts.edu}

\author[Steve Tschantz]{Steve Tschantz}
\address{Department of Mathematics\\
        Vanderbilt University\\
        Nashville, TN 32340}
\email{tschantz@math.vanderbilt.edu}

\date{\today}


\begin{abstract}
We provide conditions on the defining graph of a right-angled Coxeter group
presentation that guarantees the boundary of any $\cat(0)$ space on
which the group acts geometrically will be locally connected.
\end{abstract}

\maketitle

\centerline{0. Introduction}

\medskip

This is an edited version of our published paper \cite{MihalikRuaneTschantz}. The changes are minor, but clean up the paper quite a bit. The proof of lemma 5.8  is reworded and a figure is added showing where (in the published version) some undefined vertices of the Cayley graph are located. The proof of theorem \ref{new} is simplified by using lemma \ref{tracking} allowing us to eliminate four lemmas that appeared near the end of section 1. The statement of Lemma \ref{mainlem} is slightly changed, although the proof remains exactly the same. The revised statement of Lemma \ref{mainlem} is what is actually used in the proof of the Main Theorem (\ref{mainthm}).

In \cite{BestvinaMess}, the authors ask whether all one-ended word hyperbolic
groups have locally connected boundary.  In that paper, they relate
the existence of global cut points in the boundary to local
connectivity of the boundary.  In \cite{Bowditch}, Bowditch gives a 
correspondence between local cut points in the
boundary and splittings of the group over one-ended subgroups.  In
\cite{Swarup}, Swarup uses the work of Bowditch and others to prove that
the boundary of such a group cannot contain a global cut point, thus
proving these boundaries are locally connected.

The situation is quite different in the setting of $\cat(0)$
groups. If $G$ is a one-ended group acting geometrically on a
$\cat(0)$ space $X$, then $\partial X$ can indeed be non-locally
connected.  For example, consider the group $G=F_2\times\mathbb Z$
where $F_2$ denotes the non-abelian free group of rank 2, acting
on $X=T\times\mathbb R$ where $T$ is the tree of valence 4 (or
the Cayley graph of $F_2$ with the standard generating set).  It
is easy to see that $\partial X=\Sigma(C)$ (where $\Sigma$ denotes
the unreduced suspension), for $C$ a Cantor set and thus
$\partial X$ is not  locally connected.  Another way of viewing this group
is as an amalgamated product of two copies of $\mathbb Z\oplus\mathbb Z$
over a $\mathbb Z$ subgroup with presentation $\langle a,b,c,d | ab=ba, cd=dc, b=c\rangle$.  
This group naturally acts on a $\cat(0)$ space $X$
which is a union of planes and strips glued according to the amalgamation - in fact, 
$X$ is simply the universal cover of the presentation two complex where each square is given the metric of a unit square in the plane.   In \cite{MihalikRuane},
we prove it is easy to construct these types of examples using
amalgamated products that are ``geometric" in some sense.  In
particular, a more general version of the main theorem of that paper is given here
as Theorem 3.3 and the above example clearly satisfies the hypotheses.  

\medskip

\noindent{\bf Theorem 3.3.} {\it Suppose $A$, $B$ and $C$ are finitely generated
groups and $G=A\ast _C B$ acts geometrically on a
$CAT(0)$ space $X$.  If the following conditions are satisfied,
then $\partial X$ is not locally connected:
\begin{enumerate}
\item $[A:C]\geq 2,\ [B:C]\geq 3$.
\item There exists $s\in G-C$ with $s^n\not\in C$ for all $n\neq
0$ and $sCs^{-1}\subset C$.
\item $Cx_0$ is quasi-convex in $X$ for a basepoint $x_0$.
\end{enumerate}}

\smallskip

A general question arising naturally here is whether there is a
converse to this theorem.  In other words, is there a relationship
between the topology of the boundary and certain ``geometric"
splittings of the groups?  It is important to note here that the
boundary of a $\cat(0)$ group is not a well-defined object as is the
case for a word hyperbolic group.  Indeed, in \cite{CrokeKleiner} there is an
example of a group $G$ acting geometrically on two $\cat(0)$ spaces
with non-homeomorphic boundaries.  Nonetheless, any possible boundary
for this group will be non-locally connected because this group admits a
``geometric" splitting as in the above theorem.  

Obtaining a general converse for $\cat(0)$ groups is a very difficult
question and as of yet, there is no machinery in place to prove such
a result, thus it is worthwhile to test whether such a result holds
for a class of $\cat(0)$ groups that have some structure to work
with.  In this paper, we consider this question for right-angled
Coxeter groups. Right-angled Coxeter groups were shown to be $\cat(0)$
groups in \cite{Gromov} while general Coxeter groups were shown to be
$\cat(0)$ in \cite{Moussong}.  

In \cite{MihalikRuane}, the authors give a consequence of Theorem 3.3 above for testing whether the $\cat(0)$ Coxeter complex for a general Coxeter group has non-locally connected boundary using a presentation graph for the group.  Suppose $G$ is a Coxeter group and $S$ is a generating set for $G$ which gives rise to a Coxeter presentation $P=P(G,S)$ for $G$. The presentation can be encoded in a graph $\Gamma=\Gamma(G,S)$ called the {\it presentation graph}, whose vertex set is $S$ and edge set comes from the relations of $P$ (see  Section 1 for details).   Loosely, the consequence says (it is stated below as Theorem ~\ref{new}), if certain graph-theoretic conditions hold for $\Gamma$, then the boundary of the $\cat(0)$ Coxeter complex contains points of non-local connectivity.  The graph-theoretic conditions are the weakest possible conditions to guarantee that the Coxeter group splits as an amalgamated product as in Theorem 3.3. where the factors in the splitting are subgroups generated by subsets of the generating set.  It is Theorem ~\ref{new} for which we are trying to obtain a converse in the right-angled case.   Our first step is to show that if the graph-theoretic conditions hold, then the boundary of any $\cat(0)$ space on which the group acts will have non-locally connected boundary, not just the boundary of the $\cat(0)$ Coxeter complex.  This is the contents of Theorem ~\ref{new}.   An improved version of Theorem 3.3 is necessary to obtain this so we provide proofs of both theorems here in Section 3.    These results can be omitted if one is only interested in the main theorem.  

The statement of the main theorem is given below.  To understand the
statement of the main theorem, we briefly explain notation and terminology
here, but give formal definitions in the paper.  

If $(G,S)$ is a Coxeter system with finite generating set, then we
will denote the presentation graph by $\Gamma
(G,S)=\Gamma(S)=\Gamma$ if the $G$ and $S$ are clear.  We denote
the Cayley graph by $\Lambda(G,S)=\Lambda$.    We say a subset 
$C$ of $S$ is a {\it product separator} of $\Gamma$ if $C=A\cup B$ with $A\cap
B=\emptyset$, $\langle C\rangle = \langle A\cup B\rangle=\langle
A\rangle \oplus\langle B\rangle$, $\langle A\rangle$ and $\langle B\rangle$ are infinite and $C$ separates the graph
$\Gamma$.  A {\it virtual factor separator} (VFS) is a triple
$(C,C_1,K)$ where $C$ separates $\Gamma$, $\langle C_1\rangle$ 
has finite index in $\langle C\rangle$, $\langle K\rangle$ is
infinite and all the letters of $K$ commute with the letters of
$C_1$.  We may reduce the general problem of interest to the case when $G$ 
does not {\it visually} split as a direct product.    This means 
the graph $\Gamma$ cannot be decomposed into two non-empty disjoint
sets where each vertex of one set is joined to every vertex of
the other set and vice versa.  Thus we can assume not all the
letters of $S-C$ commute with all the letters in $C$.  And finally, a 
{\it suspended separator} is a special type of VFS where $S=C\cup\{s,t\}$ 
and $\Gamma$ (as a graph) is the suspension of $C$ with $s$ and $t$ as the suspension 
points.  In this case, $(C,C,\{s,t\})$ forms a VFS.  

\smallskip

\noindent{\bf Main Theorem.} {\it Suppose $(G,S)$ is a right-angled Coxeter system, $G$ is
one-ended, $\Gamma (S)$ contains no product separator and no VFS. Also assume that $G$ 
does not visually split as a non-trivial direct product with infinite factors. Then $G$ has locally connected boundary.}

\medskip

Note, if $G$ has more than one end, then it is trivial to determine if $G$ has locally connected
boundary or not so it is reasonable to assume $G$ is one-ended.  If $G$ splits as a
visual direct product, then $G$ has non-locally connected
boundary iff one of the factor groups does, so for the entire
paper, we assume the group $G$ does not visually split as a direct
product unless we specifically indicate that it does.    We now state Theorem ~\ref{new} which is a consequence of Theorem 3.3 above.  This result sheds light on why the other hypotheses are necessary in the Main Theorem.

\medskip

\noindent{\bf Theorem 3.2.} {\it Suppose $(G,S)$ is a Coxeter group with presentation graph
$\Gamma$.
\begin{enumerate}
\item If $\Gamma$ has a suspended separator $C$, then $G$ has locally
connected boundary if and only if $\langle C\rangle$ has locally connected boundary.
\item If  $\Gamma$ has a VFS $(C,C_1,K)$ and $C$ is not a suspended separator,
then $G$ has non-locally connected boundary.
\end{enumerate}}

To see the extent to which our Main Theorem gives a converse to Theorem 3.3, we point out the following combination of the Main Theorem and Theorem ~\ref{new}.

\medskip

\noindent{\bf Corollary.} {\it Suppose $(G,S)$ is a right-angled Coxeter system, $G$ is one-ended and $G$ does not visually split as a non-trivial direct product.  If the presentation graph $\Gamma$ contains a VFS, then $G$ has non-locally connected boundary.  If $G$ contains no product separator and no VFS, then $G$ has locally connected boundary.}  

\medskip

In Sections 4 and 5, the main technical tools for the proof are
developed.  From this point on, all our groups are right-angled since our methods do not 
work in general Coxeter groups.  Section four contains two important lemmas which allow us
to approximate ``close'' $\cat(0)$ geodesic rays by Cayley graph geodesic
rays that have a common initial subsegment (lemmas ~\ref{tracking} and ~\ref{best}).  
In general $\cat(0)$ group theory, one cannot hope to use Cayley graph geometry to study the geometry of geodesics in the $\cat(0)$ space since the natural quasi-isometry between them is badly behaved.   In this paper, we prove that if the group is a right-angled Coxeter group, we can gain enough control over this quasi-isometry to use the intuition coming from the Cayley graph geometry inside the $\cat(0)$ space.  This is truly the main technical difficulty in the paper.

Section 5 contains the construction of a {\it filter} for a pair of Cayley
graph geodesics needed to prove the main theorem.  In doing so, we build a
collection of graphs in the plane whose edges are labeled by elements of the generating set
$S$.  There is an obvious map from any such graph into the Cayley graph $\Lambda$. A
final version of these graphs will be 1-ended and the natural map from this
graph into $\Lambda$ will be {\it proper}.  If the two geodesics we start with have a long
common initial piece, we will show that this filter maps over to $X$ under the natural quasi-isometry
in a controlled manner.  

In Section 6, we give the proof of the main theorem.   Our goal is to prove local connectivity of the boundary so we start by considering two ``close'' geodesic rays (equivalently, two close boundary points) in the $\cat(0)$ space.  From these, we obtain two approximating Cayley graph geodesic which share a long common subpiece.  Using Section 5, we build a filter between these two rays which has a nice property.  In particular, the graph-theoretic hypotheses allow us to conclude there is a global bound on the length of  a {\it factor path} in the filter (see lemma ~\ref{mainlem}).  A factor path is an edge path in the Cayley graph that lies in a product subgroup (see definition ~\ref{factor}) - these are exactly the paths that behave badly under quasi-isometry.  Carrying this information back to $X$ via the natural quasi-isometry allows us to construct a ``small'' connected set in the boundary of the space containing the two points and hence show the $\cat(0)$ space has locally connected boundary.    
\newpage
\section {Coxeter group preliminaries}

In this section, we prove several technical facts about Coxeter groups.  Many of these facts can be found scattered about in the vast literature on Coxeter groups, however we state them here with our notation and with our particular viewpoint in mind.  Even though Theorem ~\ref{mainthm} only holds for right-angled Coxeter goups, we prove the results in this section assuming $G$ is a general Coxeter group.  We will need the more general results in section 3 of this paper.  

There are two types of results in this section.  The results up through lemma ~\ref{relator} give geometric information about the Cayley graph of a Coxeter group.  The most important lemmas for understanding the proof of the main theorem are theorem ~\ref{dc} and lemmas ~\ref{convex}, ~\ref{backward}, ~\ref{mike}.  We need these geometric results in sections 4 and 5 where we must compare the geometry of the $\cat(0)$ space $X$ with the combinatorial properties of the Cayley graph $\Lambda$.  The remaining results culminate with lemma ~\ref{finiteind}.  This final lemma says if $H$ is a {\it special} subgroup of $G$ of finite index, then $G$ splits as a direct product in a way that you can {\it see} in the presentation graph of $G$.  This lemma is used to prove lemma ~\ref{graph} which is necessary to begin building a filter for a pair of geodesic rays. 

A {\it Coxeter group} is a group $G$ with generating set
$S=\{s_{1},\ldots s_{n}\}$
having a {\it Coxeter presentation} of the form $P(G,S) =\langle S\
\vert\
(s_{i}s_{j})^{m_{ij}}=1\rangle $ for $1\le i,j\le n$ with $m_{ij}=1$ if and
only if
$i=j$, otherwise $2\le  m_{ij}\le\infty$.  If $m_{ij}=\infty$, then
$s_{i},s_{j}$ have no relation.  The pair $(G,S)$ is called a {\it Coxeter
System}. If $m_{ij}=2 \ {\text or} \ \infty$ for all relations in the
presentation, $G$ is called {\it right-angled}.

The {\it presentation graph} of $(G,S)$ is the graph $\Gamma (G,S)$ with
vertex set $S$ and an edge labeled $m$ between vertices $s$ and $t$ if
$(st)^{m}=1$ is a relation of $P$.  For right-angled groups, we omit
the $m$ labels, thus two generators commute if and only if there is
an edge between them in the presentation graph.

Assume $(G,S)$ is a Coxeter system with $S$ a finite set.  Let
$\Lambda (G,S)=\Lambda$ be
the Cayley graph of $G$ with respect to the generating set $S$.
We assume each edge of $\Lambda$ is labeled by an element
of $S$. Each $s\in S$ is an element of order two and has a fixed
point set in $\Lambda$ denoted by ${\text Fix}(s)$.  This
set consists  of the midpoints of edges labeled by $s$ with each
endpoint in the subgroup of $G$ generated by
the letters of $S$ that are adjacent to $s$ in $\Gamma (G,S))$.

An {\it edge path} in $\Lambda$ can be denoted as $(e_{1},\ldots
e_{n})$ where $e_{i}\in S$ for all $i$ as long as the initial vertex
of the path is specified. In many of our arguments the initial vertex is
not important as the argument is valid for such a path at any vertex.
As an example, whether or not a path is geodesic does not depend on
the initial vertex.
In these cases we may supress initial vertices.

Perhaps the most useful combinatorial and geometric fact about Coxeter groups is the
following which is referred to as the {\it Deletion Condition} in this paper.  It is used in this paper 
in any result that involves the geometry of the Cayley graph $\Lambda$.

\begin{thm}\label{dc} {\bf (The Deletion Condition)} If $w=w_{1}\ldots w_{n}$ is
non-geodesic
in $\Lambda$ then there are $0<i<j\leq n$ such that $$w=w_{1}\ldots
w_{i-1}w_{i+1}\ldots w_{j-1}w_{j+1}\ldots w_{n}.$$
\end{thm}

\noindent{\bf Remark.} In this case, we say that $w_i$ {\it deletes}
with $w_j$.
Observe that if $w_1\ldots w_{n-1}$ is geodesic, but $w_1\ldots
w_{n-1}w_n$ is not
geodesic, then we know $w_n$ must delete with one of the previous
letters.  If the group
is right-angled, then $w_n$ deletes with $w_i$ for some $1\le i< n$.
In fact, $w_i=w_n$
and $i$ is the largest integer between 1 and $n-1$ with $w_i=w_n$.
Also, $w_n$ commutes
with $w_j$ for all
$i<j\le n-1$.

\begin{defn}[Special subgroup] Suppose $(G,S)$ is a Coxeter system and let $V\subset
S$. Then the subgroup of $G$ generated by $V$ is called a {\it special
subgroup} of $G$.
\end{defn}

The next result shows why the special subgroups are indeed special.

\begin{thm} {\bf (Special Subgroup Theorem)} If $(G,S)$ is a Coxeter
system and $V\subset S$ then the subgroup $\langle V\rangle$ is
Coxeter with presentation $\langle V\vert R\rangle $ where $R$ is
the set of all relations of $P(G,S)$ involving only letters from
$V$.
\end{thm}

The next two lemmas are part of the standard theory of Coxeter groups
and proofs can be found in \cite{Humphreys} for example. The first of these 
lemmas says that special subgroups are convex in $\Lambda(G,S)$.  

\begin{lemma}\label{convex} If $U\subset S$ and $u\in \langle U\rangle $ is written
as a geodesic word in $U$ then this word is geodesic in $S$.
Furthermore every geodesic for $u$ can only use letters from $U$.
\end{lemma}

If $\alpha$ is a path in $\Lambda$ then define $\bar \alpha \in G$
to be the product of the edge labels of $\alpha$.  The following lemma says that if we add one 
more edge to a geodesic path of length $n$ in $\Lambda$, then the new path cannot also be of length $n$.  Also, the set of generators that make this new path have length $n-1$ only generate a finite special subgroup.  This will be very important in lemmas ~\ref{better} and ~\ref{finite}. 

\begin{lemma}\label{backward} If $\alpha$ is geodesic in $\Lambda$ then the set of
labels of edges $e$ such that $( \alpha ,e )$ is not geodesic
generates a finite subgroup of $G$, denoted $B(\bar \alpha)$.
Furthermore, either $l(\alpha e)=l(\alpha)-1$ or $l(\alpha
e)=l(\alpha)+1$
\end{lemma}

The next three lemmas are geometric consequences of the convexity of the special subgroups which will be useful in the proof of the main theorem.  The first of these follows directly from the Deletion Condition and lemma ~\ref{convex}. It allows us to construct geodesics in $\Lambda$ by combining geodesics from special subgroups with non-overlapping generating sets. 

\begin{lemma}\label{comb} Suppose $(G,S)$ is a Coxeter
system, $U$ and $V$ are subsets of $S$ with $U\cap V=\emptyset$.  If $\alpha$ is a
geodesic in the letters of $U$ and $\beta$ is a geodesic in the letters of $V$, then
$(\alpha,\beta)$
is a geodesic in $\Lambda(G,S)$.
\end{lemma}

The next lemma is essentially a more technical way of saying the special subgroups are convex.  Specifically, if you start with a word of the form $sut$ where $s,t\in S$ with $t,u$ both in a particular special subgroup (the one generated by the letters of $u$).  Then if the edge path $sut$ ends in the subgroup, $s$ must also be in that subgroup since $s$ must equal $t$ and this generator commutes with the letters of the word $u$.  This form of the convexity is used in many of the results of this section and in lemmas ~\ref{better} and ~\ref{sep}.

\medskip

\begin{lemma}\label{mike}Suppose $U\subset S$. If $s\in S-U$, $t\in S$ and for
some $u\in
\langle U\rangle$, $sut\in \langle U\rangle $, then
\begin{enumerate}
\item $t=s$.
\item $sus=u$ and
\item If $u=u_{1}\ldots u_{n}$ is geodesic, (so $u_{i}\in U$),
then $su_{i}s=u_{i}$ for all $i$.
\end{enumerate}
\end{lemma}

\begin{proof} As $sut\in \langle U\rangle $, $s\in \langle U\cup
\lbrace t\rbrace \rangle $. By lemma ~\ref{convex}, $t=s$ and (1) is finished.

Say $u=u_{1}\ldots u_{n}$ is geodesic, where $u_{i}\in U$. By
lemma ~\ref{convex}, this is geodesic in the entire Coxeter system.
Say $sus=v_{1}\ldots v_{m}$, $m$ minimal and $v_{i}\in U$.
The words $su_{1}\ldots u_{n}$ and $v_{1}\ldots v_{m}s$ are geodesics (for the same element of $G$) by lemma \ref{comb} and so $m=n$.
Now $su_{1}\ldots u_{n}s$ is not geodesic. The second $s$ cannot
delete with a $u_{i}$ as $s\not\in \langle U\rangle $. So $s$ commutes
with $u$, finishing (2).

To prove (3), we proceed  by induction on $n$. If $n=1$, (3) is clear
from (2).
Assume (3) is true for $k\leq n-1$ and $n\geq 2$. The path
$u_{1}\ldots u_{n}s$ is geodesic (as $s\not\in U$), but $u_{1}\ldots
u_{n}su_{n}$ ($=su_{1}\ldots u_{n-1}$) is not geodesic. Clearly the
last $u_{n}$ and $s$ do not delete as $s\not=u_{n}$. Say the last
$u_{n}$ and $u_{i}$ delete. This implies (by part (2)) that
$s$ commutes with both $u_{1}\ldots u_{n-1}$ and $u_{1},\ldots
u_{i-1}u_{i+1}\ldots u_{n}$. By induction, $s$ commutes
with $u_{i}$ for all $i$.
\end{proof}

It is an elementary exercise using the deletion condition to show that because of convexity of the special subgroups, we also have the notion of a unique projection onto a special subgroup (or one of it's cosets).  

\begin{lemma}\label{proj} Suppose $\ast$ and $v$ are vertices of $\Lambda$ and
$T\subset S$. Then there is a unique vertex $w$ of $v\langle T\rangle$ closest to $\ast$.
Furthermore, if $\alpha$ is a geodesic from $\ast$ to $w$ and $\beta$
is a geodesic in the letters of $T$, then $(\alpha,\beta)$ is a geodesic in
$\Lambda$.
\end{lemma}

Lemma ~\ref{proj} and the Deletion Condition now imply:

\begin{lemma}\label{relator} Suppose $\beta$ is a geodesic from $x$ to $y$ in
$\Lambda$ and $( \beta , s)$ is not geodesic. If $t$ is the
last letter of $\beta$ then $s$ is related to $t$ and the edge path
labeled by half of this relation, at $y$, is geodesic towards $x$.
\end{lemma}

\begin{defn}[2-link] \label{link}For $A\subset S$, we define $lk^2(A)=\cap_{a\in A}lk^{2}(a)$ where $lk^{2}(a)$
means all $s\in S$ that are connected to $a$ in the graph
$\Gamma(G,S)$ with an edge labeled 2.  In particular, for $s\in S$, we often write $lk^2(s)$ instead of $lk^2(\{s\})$.   
\end{defn}

We make the following observations concerning the right-angled Coxeter groups.  These statements will be used in the last three sections of the paper when we restrict to the right-angled case.

\medskip

\noindent{\bf Remark.} Let $(G,S)$ be right-angled and $H\subset G$. Let $C(H)$ denote the centralizer of $H$ in $G$. 
\begin{enumerate}
\item   For $(G,S)$ right-angled with $A\subset S$,
$lk^2(A)=lk(A)$ where $lk(A)$ is the usual link of $A$ in the graph $\Gamma (G,S)$. 
\item For $s\in S$, $C(s)\slash\langle s\rangle=\langle lk(s)\rangle$.  This follows
from lemma ~\ref{mike} and the Deletion Condition.  In particular, $\langle
lk(s)\rangle$ is an index two subgroup of $C(s)$ in this case.
\end{enumerate}

\begin{lemma}\label{ext}  If $V$ does not visually split
non-trivially as a direct product, then for any geodesic $\alpha$ in $V$ and $v\in V$,
$\alpha$ can be extended to a geodesic ending with $v$.
\end{lemma}

\begin{proof} If not, choose $\alpha$ a geodesic such that
$Y_{\alpha}=\{y\in S|\
\alpha$ cannot be extended to a geodesic ending with $y\}$ is as
large as
possible.  By hypothesis, $Y_{\alpha}\ne\emptyset$ and by lemma ~\ref{backward},
$Y_{\alpha}\ne
S$.  By the maximality of $Y_{\alpha}$, $\alpha$ can be extended to
an infinite
geodesic $(\alpha,\gamma)$ such that the letters of $S-Y_{\alpha}$
occur
infinitely often in $\gamma$.  We show $Y_{\alpha}$ commutes with
$S-Y_{\alpha}$ to
obtain a contradiction.

Assume $\alpha=(a_1,\ldots,a_n)$ and $\gamma=(b_1,b_2,\ldots)$.  Let
$\beta_i=(\alpha,b_1,b_2,\ldots,b_i)$.  For $y\in Y_{\alpha}$,
$(\beta_i,y)$ is not
geodesic.  This last $y$ cannot delete with a $b_j$ as $y\notin
S-Y_{\alpha}$ and so
must delete with some $a_k$.  Choose $i< j$ such that all letters of
$S-Y_{\alpha}$
occur in $\{b_{i+1},\ldots, b_j\}$ and the $y$ ending $(\beta_i,y)$
and $(\beta_j,y)$
both delete with $a_k$.  Then $a_k\ldots a_nb_1\ldots
b_i=a_{k+1}\ldots a_nb_1\ldots
b_iy$ and $a_k\ldots a_nb_1\ldots b_jy=a_{k+1}\ldots a_nb_1\ldots
b_j$.  So
$a_{k+1}\ldots a_nb_1\ldots b_iyb_{i+1}\ldots b_jy=a_{k+1}\ldots
a_nb_1\ldots b_j$ and
$yb_{i+1}\ldots b_jy=b_{i+1}\ldots b_j$.  By lemma ~\ref{mike}, $y$ commutes
with
$S-Y_{\alpha}$ as needed. 
\end{proof}

\begin{lemma}\label{finiteind} {\bf (Finite Index Lemma)} Suppose $(G,S)$ is a Coxeter
system,
$G$ is infinite, $U\subset S$, $\langle U\rangle $ and has finite
index
in $G$ and $F\subset U$ is the maximal set such that $\langle
F\rangle$
is finite and $\langle U\rangle$ visually splits as $\langle F\rangle \oplus \langle
U-F\rangle $ then $\langle (S-U)\cup F\rangle $ is finite and $G$
splits as $\langle (S-U)\cup F\rangle \oplus \langle U-F\rangle $.
\end{lemma}

\begin{proof} Clearly $\langle U-F\rangle$ has finite index in $V$.
Hence there is a
bound $K$ such that every vertex of $\Lambda(G,S)$ is within $K$ of
$\langle
U-F\rangle$.  If $\langle(S-U)\cup F\rangle$ were not finite, choose
a geodesic $\alpha$
of length $K+1$ in the letters $(S-U)\cup F$, with initial point the
identity.  The end
point $y$ of this geodesic is within $K$ of a vertex $x\in\langle
U-F\rangle$.  Let
$\beta$ be a geodesic in the letters of $U-F$ from $x$ to the
identity.  By the
lemma ~\ref{comb}, $(\beta, \alpha)$ is a geodesic from $x$ to
$y$, but
$d(x,y)\le K$, giving the desired contradiction.  So
$\langle(S-U)\cup F\rangle$ is
finite.

Maximally decompose $\langle U-F\rangle$ as $\oplus_{i=1}^n\langle
W_i\rangle$ with
$W_i\ne\emptyset$, $\cup_{i=1}^n W_i=U-F$, and $W_i\cap
W_j=\emptyset$ for $i\ne j$.
By the lemma ~\ref{ext}, there is a geodesic $\alpha_i$ with
letters in $W_i$
that begins at $1\in\Lambda$ such that each letter of $W_i$ occurs
infinitely often in
$\alpha_i$.  Let $\beta$ be an edge path at $1$ such that each
element of $(S-U)\cup
F$ labels exactly one edge of $\beta$ and $\beta$ has length ${\text
Card}((S-U)\cup
F)$.  By lemma ~\ref{convex}, $\beta$ is geodesic.  Say $\alpha_i=(a_{i1},
a_{i2},\ldots)$ and
$\beta=(b_1,\ldots,b_n)$.  Choose $k$ such that all letters of $W_i$
occur at least
$n+K+1$ times in $\alpha_{ik}=(a_{i1},\ldots,a_{ik})$.  Let $x$ be
the end point of
$\beta$ and $y$ the end point of $\alpha_{ik}$.  Let $\delta$ be a
geodesic from $y$ to
$x\langle U-F\rangle$ of length $\le K$.  Say $\delta$ ends at $z$.  Let
$\gamma=(c_1,\ldots,
c_m)$ be a geodesic from $x$ to $z$ with $c_j\in (U-F)$ for all $j$.
By the
lemma ~\ref{comb} $(\alpha_i^{-1},\beta)$ is geodesic and by lemma
~\ref{proj},
$(\delta,\gamma^{-1})$ is also geodesic.  Hence
$(\alpha_{ik}^{-1},\beta,c_1)$ is not
geodesic.

Say $c_1$ deletes wtih $a_{ij}$.  Since $b_i\notin U-F$ for all $i$,
we must have $c_1$
commuting with $b_j$ for all $j$ by lemma ~\ref{mike}.  From the Deletion
Condition we see that
if $(u_1\ldots, u_n)$ and $(v_1,\ldots,v_n)$ are geodesics with the
same endpoints,
then $\{u_1,\ldots,u_n\}=\{v_1,\ldots,v_n\}$, thus we obtain
$\{c_1,\ldots,c_m\}\subset
U-F$.  If we replace $\alpha_i$ with a geodesic in $U-F$ from $c_1$
to $y$, we see that
$c_2$ commutes with each $b_j$.  Continuing, each $b_j$ commutes with
each $c_q$.  Now
$c_1\ldots c_m=a{i1}\ldots a_{ik}\delta\beta^{-1}$.  At most $n+K$
deletions can occur
in the right hand side before a geodesic is obtained.  Hence all
letters of $U-F$
appear in the resulting reduced word.  So all letters of $U-F$ appear
in
$\{c_1\ldots,c_m\}$ forcing the letters of $U-F$ to commute with the
letters
of $(S-U)\cup F$ as needed. 
\end{proof}

\section{$\cat(0)$ spaces and their boundaries}

In this section we give definitions and basic properties of $\cat(0)$
spaces, boundaries and isometries as well as some known facts we will
need in the proof of the main result.

Let $(X,d)$ be a metric space.  Then $X$ is {\it proper} if closed metric
balls are compact.  A ({\it unit speed}){\it geodesic} from $x$ to $y$ for
$x,y\in X$ is a map $c:[0,D]\to X$ such that $c(0)=x,\ c(D)=y$ and
$d(c(t),c(t'))=|t-t'|$ for all $t,t'\in [0,D]$.  If
$I\subseteq\mathbb R$ then a map $c:I\to X$ parametrizes its image
{\it proportional to arclength} if there exists a
constant $\lambda$ such that $d(c(t),c(t'))=\lambda|t-t'|$ for all
$t,t'\in I$.  Lastly, $(X,d)$ is a called a {\it geodesic metric
space} if every pair of points are joined by a geodesic.

\begin{defn}[$\cat(0)$]  Let $(X,d)$ be a proper complete geodesic
metric space.
If $\vartriangle abc$ is a geodesic triangle in $X$, then we consider
$\vartriangle\overline a\overline b\overline c$ in $\mathbb E^2$, a
triangle with the same side lengths, and call this a {\it comparison
triangle}.  Then we say $X$ satisfies the $\cat(0)$ {\it inequality}
if given $\vartriangle abc$ in $X$, then for any comparison triangle and
any two points $p,q$ on $\vartriangle abc$, the corresponding points
$\overline p,\overline q$ on the comparison triangle satisfy
$$d(p,q)\leq d(\overline p,\overline q)$$
\end{defn}

If $(X,d)$ is a $\cat(0)$ space, then the following basic properties
hold:
\begin{enumerate}
\item  The distance function $d\colon X\times X\to\mathbb R$ is convex.
\item  $X$ has unique geodesic segments between points.
\item  $X$ is contractible.
\end{enumerate}
For details, see \cite{BridsonHaefliger}.

Let $(X,d)$ be a proper $\cat(0)$ space. First, define the boundary,
$\partial X$ as a set as follows:

\begin{defn}[Asymptotic] Two geodesic rays $c,c^\prime\colon [0,\infty)\to
X$ are said to be {\it asymptotic} if there exists a constant $K$ such
that $d(c(t),c^\prime (t))\leq K, \forall t>0$ - this is an
equivalence relation.  The boundary of X, denoted $\partial X$, is then the
set of equivalence classes of geodesic rays.  The union $X\cup\partial X$
will be denoted $\overline X$.  The equivalence class of a ray $c$ is
denoted by $c(\infty)$.
\end{defn}

There is a natural neighborhood basis for a point in $\partial
X$.  Let $c$ be a geodesic ray emanating from $x_0$ and $r>0,\
\epsilon >0$. Also, let $S(x_0,r)$ denote the sphere of radius $r$ centered
at $x_0$ with $p_r:X\to S(x_0,r)$ denoting projection.  Define
$$U(c,r,\epsilon)=\{ x\in\overline{X} \vert d(x,x_0)>r,\
d(p_r(x),c(r))< \epsilon\}$$
This consists of all points in $\overline{X}$ such that when projected
back to $S(x_0,r)$, this projection is not more than $\epsilon$ away
from the intersection of that sphere with $c$.  These sets along with the
metric balls about $x_0$ form a basis for the {\it cone topology}.
The set $\partial X$ with the cone topology is often called the
{\it visual boundary}.  As one expects, the visual boundary of
$\mathbb R^n$ is $S^{n-1}$ as is the visual boundary of $\mathbb H^n$.

\begin{defn} Let $\gamma$ be an isometry of the metric space
$X$.  The displacement function $d_{\gamma}\colon
X\to\mathbb\mathbb R_+$ is defined by $d_{\gamma}(x)=d(\gamma\cdot x,x)$.  The
translation length of $\gamma$ is the number $|\gamma|=
inf\{d_{\gamma}(x): x\in X\}$.  The set of points where $\gamma$
attains this infimum will be denoted $\text{Min}(\gamma)$.  An
isometry $\gamma$ is called {\it semi-simple} if $\text{Min}(\gamma)$ is
non-empty.
\end{defn}

We summarize some basic properties about this $\text{Min}(\gamma)$ in
the following proposition.

\begin{proposition}  Let $X$ be a metric space and $\gamma$ an
isometry of $X$.
\begin{enumerate}
\item $\text{Min}(\gamma)$ is $\gamma$-invariant.
\item If $\alpha$ is another isometry of $X$, then
$|\gamma|=|\alpha\gamma\alpha^{-1}|$, and
$\text{Min}(\alpha\gamma\alpha^{-1})=\alpha\cdot \text{Min}(\gamma)$;
in particular,
if $\alpha$ commutes with $\gamma$, then it leaves
$\text{Min}(\gamma)$
invariant.
\item If $X$ is $\cat(0)$, then the displacement function $d_{\gamma}$
is convex: hence $\text{Min}(\gamma)$ is a closed convex subset of
$X$.
\end{enumerate}
\end{proposition}

\begin{defn} Let $X$ be a metric space.  An isometry
$\gamma$ of $X$ is called
\begin{enumerate}
\item elliptic if $\gamma$ has a fixed point - i.e $|\gamma|=0$ and
$\text{Min}(\gamma)$ is non-empty.
\item hyperbolic if $d_{\gamma}$ attains a strictly positive infimum.
\item parabolic if $d_{\gamma}$ does not attain its infimum, in other
words if $\text{Min}(\gamma)$ is empty.
\end{enumerate}
\end{defn}

It is clear that an isometry is semi-simple if and only if it is
elliptic or hyperbolic.  If two isometries are conjugate in
$\text{Isom}(X)$, then they are in the same class.

If a group $\Gamma$ acts geometrically on a $\cat(0)$ space $X$, then
the elements of $\Gamma$ act as semi-simple isometries because of the
cocompactness of the action.

The next theorem implies that the centralizer of an element $\gamma$
in a $\cat(0)$
group $\Gamma$ is again a $\cat(0)$ group because it acts
geometrically on $\cat(0)$
subspace ${\text Min}(\gamma)$.  A proof of this can be found in
\cite{Ruane}.

\begin{thm}\label{kim}  Suppose $\Gamma$ acts geometrically on the
$\cat(0)$ space $X$ and suppose $\gamma\in\Gamma$.  Then
$C_{\Gamma}(\gamma)$ acts geometrically on the  $\cat(0)$ subset
$Min(\gamma)$ of $X$.
\end{thm}

\noindent{Remarks.}  We make the following remarks which will be needed in our setting of right-angled Coxeter groups.
\begin{enumerate}
\item If $(G,S)$ is a right-angled Coxeter group
acting geometrically on the $\cat(0)$ space $X$, then each $s\in S$ is an
elliptic isometry with ${\text Min}(s)={\text Fix}(s)$ a closed, convex
subset of $X$ on which the centralizer of $s$ acts geometrically.
\item Since $C(s)$ acts geometrically on ${\text
Fix}(s)$ and $\langle lk(s)\rangle$ is an index two subgroup of $C(s)$ (see remarks
following definition ~\ref{link}) thus $\langle lk(s)\rangle$ is a $\cat(0)$
subgroup of $G$.
\end{enumerate}

We will need the following result concerning $\cat(0)$ products.
Proofs of these can be found in \cite{BridsonHaefliger}.

\begin{thm}\label{product} Suppose $X$ and $Y$ are $\cat(0)$ spaces.Then $X\times Y$ with the product
metric is also a $\cat(0)$ space.  Furthermore, $\partial(X\times
Y)$ is homeomorphic to the spherical join of $\partial X$ and
$\partial Y$.  In particular, if $Y=\mathbb R$, then
$\partial(X\times\mathbb R)$ is homeomorphic to the suspension of
$\partial X$.
\end{thm}

The following theorem is due to Milnor in \cite{Milnor} and will be used
throughout the paper to carry information between the Cayley graph of
$(G,S)$ and the $\cat(0)$ space $X$.

\begin{thm}\label{milnor} If a group $G$ with finite generating set $S$ acts
geometrically on a proper, geodesic metric space X, then the Cayley
graph of $G$ with
respect to $S$ is quasi-isometric to $X$ under the map $g\mapsto
g\cdot x_0$ where
$x_0$ is a fixed basepoint in $X$.
\end{thm}

\section{Local Connectivity}

This section contains the proof of Theorem ~\ref{new} below.  The
material in this section is not needed in the proof of the main
theorem.  We say a $\cat(0)$ group $G$ has {\it (non-)locally
connected boundary} if for every $\cat(0)$ space $X$ on which $G$
acts geometrically, $\partial X$ is (non-)locally connected.
Clearly, if $G$ is infinite ended, then $G$ has non-locally
connected boundary, thus we only need to consider Coxeter groups
$(G,S)$ for which the presentation graph $\Gamma$ is connected.

\begin{defn}[VFS]\label{VFS} Suppose $(G,S)$ is a Coxeter system.  Let $\Gamma$
denote the presentation graph for $(G,S)$.  A {\it virtual factor
separator} (VFS) for $\Gamma$ is a triple $(C,C_{1},K)$ where
$C_{1}\subset C$ are full subgraphs of $\Gamma$, $C$ separates
$\Gamma$, the group $\langle C_{1}\rangle $ has finite index in
$\langle C\rangle $, $K\subset lk^{2}(C_{1})$, and $\langle K\rangle
$ is infinite. If $S=C\cup \{s,t\}$, $s$ and $t$ do not span
an edge of $\Gamma$ and  $s,t\in lk^2(C)$, then $C$ is called a
{\it suspended separator} of $\Gamma$. Note that if $C$ is a
suspended separator of $\Gamma$, then $G=\langle C\rangle \oplus
\langle s,t\rangle$ and $(C,C,\{s,t\})$ is a VFS.
\end{defn}

\begin{thm} \label{new}Suppose $(G,S)$ is a Coxeter group with presentation graph
$\Gamma$.
\begin{enumerate}
\item If $\Gamma$ has a suspended separator $C$, then $G$ has locally
connected boundary if and only if $\langle C\rangle$ has locally connected boundary.
\item If $\Gamma$ has a VFS $(C,C_1,K)$ and $C$ is not a suspended separator,
then $G$ has non-locally connected boundary.
\end{enumerate}
\end{thm}

The proof of Theorem ~\ref{new} will be an application of the following
Theorem from \cite{MihalikRuane}.  The theorem cannot be applied
directly because of condition three below, but lemma \ref{tracking} of the next section implies condition three holds for special subgroups of Coxeter groups.

\begin{thm}\label{old} Suppose $A$, $B$ and $C$ are finitely generated
groups and $G=A\ast _C B$ acts geometrically on a
$CAT(0)$ space $X$.  If the following conditions are satisfied,
then $\partial X$ is not locally connected:
\begin{enumerate}
\item $[A:C]\geq 2,\ [B:C]\geq 3$.
\item There exists $s\in G-C$ with $s^n\not\in C$ for all $n\neq
0$ and $sCs^{-1}\subset C$.
\item $Cx_0$ is quasi-convex in $X$ for a basepoint $x_0$.
\end{enumerate}
\end{thm}

\noindent{\bf Remark.} If $C_{1}$ is a subgroup of
finite index in $C$, and condition (2) is replaced by ($2'$) below then the proof remains
unchanged.

\medskip

\noindent ($2'$) There exists $s\in G-C_{1}$ with $s^n\not\in C$ for
all $n\neq 0$ and
$sC_{1}s^{-1}\subset C_{1}$.

\medskip

\noindent {\it Proof of Theorem ~\ref{new}}

If $C$ is a suspended separator, then $G$ splits as $\langle
C\rangle \oplus \langle \{s,t\}\rangle$. In this case, the
element $(st)$ is virtually central -meaning, it is central in a
subgroup of finite index in $G$.  In fact, the subgroup generated
by $C$ and the element $(st)$ is a subgroup of finite index in
$G$ in which $st$ is central.  Theorem 3.4 of  \cite{Ruane}
implies that if $X$ is any $\cat(0)$ space on which $G$ acts
geometrically, then the following are true.  A subset $Z$ of $X$
is {\it quasi-dense} if  there exists a constant $K\ge 0$ such
that every $x\in X$ is within $K$ of an element of $Z$.
\begin{enumerate}
\item There exists a quasi-dense closed, convex subset $Z$ of $X$such that
$\partial Z=\partial X$ (we know $\partial Z\subset\partial X$
since $Z$ convex in $X$. Thus for boundary considerations, we
only need to know about the space $Z$.
\item $Z$ splits isometrically as $Y\times\mathbb R$ where $Y$ is a closed,
convex subset of $X$ and $Y$ admits a geometric group action by
the group $\langle C\rangle$.
\end{enumerate}

Since $Z$ splits as a product, we know $\partial Z=\partial
X\equiv\Sigma(\partial Y)$ by Theorem ~\ref{product}. Thus $\partial X$ is
locally connected if and only if $\partial Y$ is locally
connected.

Notice that if $Y$ is a $\cat(0)$ space on which $\langle
C\rangle$ acts geometrically, then we can construct the $\cat(0)$
space $X=Y\times\mathbb R$ and a geometric action of $G$ on $X$
via the product action.  Simply let $\langle C\rangle$ act on $Y$
as given and let $\langle\{s,t\}\rangle$ act by the infinite
dihedral action on the $\mathbb R$ factor.  This finishes part 1
of the theorem.

Suppose $(C,C_1,K)$ is a VFS for $\Gamma$.  Then $G=\langle
A\rangle \ast_{\langle C\rangle} \langle B\rangle$ where $A\cup
B=S$ and $A\cap B=C$ since $C$ separates $\Gamma$.  The proof of
this part of the theorem is a direct application of Theorem ~\ref{old} and lemma \ref{tracking}. 
$\square$

\section{Two important lemmas}

For the remainder of the paper, we will only be dealing with right-angled Coxeter groups.  

In this section, we prove two important lemmas about a right-angled
Coxeter group $G$ acting on a $\cat(0)$ space $X$ that will allow us to
transfer combinatorial information about the Cayley graph of $G$ to the
space $X$.  The first is lemma ~\ref{tracking} which allows us to approximate
$\cat(0)$ geodesic rays in $X$ with Cayley graph rays, and also implies that special subgroups of $G$ are quasi-convex in $X$.  The second is lemma ~\ref{best} which
says that if two $\cat(0)$ geodesic rays are sufficiently close, then the
Cayley graph approximating rays can be chosen to have a long common
subpeice.  The method of proof used here does not extend to give the
corresponding results for general Coxeter groups.

Let $\Gamma (G,S)$ denote the presentation graph for $G$ with respect to the
generating set $S$ and $\Lambda (G,S)=\Lambda$ denote the Cayley graph
of $G$ with respect to the generating set $S$.  We know each $s\in S$
is order 2 and has a closed, convex, fixed point set in $X$ denoted by
${\text Fix} (s)$.  With a slight abuse of notation, we also denote
the fixed point set of $s$ in $\Lambda$ by ${\text Fix}(s)$.  Recall
this consists of the midpoints of edges labeled by $s$ with each
endpoint in $\langle lk(s)\rangle$.

Let $x_0$ be a basepoint in $X$.  We identify a quasi-isometric
copy of the Cayley graph $\Lambda$ of $G$ inside $X$ with the
orbit of $x_0$ under the action of $G$ using Theorem ~\ref{milnor}.  We
will use $\cat(0)$ geodesics between adjacent vertices in the
Cayley graph.  With this identification, any geodesic $\alpha$ in
$\Lambda$ is assigned a piecewise $\cat(0)$ geodesic in $X$ which
we denote by $\alpha_{X}$.

We stop to make an important observation about geodesics in
$\Lambda$.

\medskip

\noindent{\bf Remark.} Suppose $\alpha=(e_1,e_2,\ldots, e_n)$ with
each $e_i\in S$ is an edge path in $\Lambda$.  The path $\alpha$ crosses a set of fixed
point sets or {\it walls} as the letters of $\alpha$ are traversed one by one.  These
walls are (in order):
\newline
${\text Fix}(e_1),{\text Fix}(e_1e_2e_1),\ldots ,{\text
Fix}(e_1e_2\cdots
e_{n-1}e_ne_{n-1}\cdots e_2e_1)$.

In other words, at the $i$th step, $\alpha$ crosses the wall
$e_1\cdots e_{n-1}{\text Fix}(e_i)$.  We refer to this wall as an
$e_i$-wall.  Notice that $\alpha$ can have the same letter $e_i$
occuring twice which gives (possibly) two different $e_i$-walls
that $\alpha$ crosses.  With this in mind, the Deletion Condition
can then be interpreted geometrically as:  the edge path $\alpha$
is geodesic in $\Lambda$ if and only if it crosses each wall at most once. 

\begin{defn}[Tracking] Suppose $r:[a,b]\to X$ is a geodesic segment in
$X$ with $r(a)=x$ and $r(b)=y$ and suppose $\delta >0$.  We say
the Cayley graph geodesic $\alpha$, $\delta$-tracks $r$ (or more
precisely, the image of $r$) if every point of $\alpha_{X}$ is
within $\delta$ of a point of the image of $r$ and the endpoints
of $r$ and $\alpha_{X}$ are within $\delta$ of each other.
\end{defn}

\begin{lemma}\label{tracking} There exists a $\delta >0$ such that for any
geodesic ray $r:[0,\infty)\to X$ based at $x_0$, there exists a
geodesic ray $\alpha_r$ in $\Lambda (G,S)$ that $\delta$-tracks $r$. Furthermore, for $A\subset S$, the special subgroup $\langle A\rangle$ is quasi-convex in $X$. 
\end{lemma}

\begin{proof} First assume $r:[0,d]\to X$ is a unit speed
parametrization of the geodesic segment in $X$ between $x_0$ and
$g\cdot x_0$ for some $g\in G$.  Since $G$ acts cocompactly on $X$,
there exists a $K>0$ such
that any point $x\in X$ is within $K$ of an orbit point.  For each
$i=1,2,\ldots D=\lfloor d \rfloor$ (i.e. $D$ is the largest integer
less than or equal to d) choose an element $g_i\in G$ such that
$d(r(i), g_i\cdot x_0)<K$.  Choose Cayley graph geodesics $\alpha_0$
from
$1_G$ to $g_1$, $\alpha_1$ from $g_1$ to $g_2$,$\ldots$,$\alpha_D$
from
$g_D$ to $g$.  Call the piecewise geodesic $\alpha$.  In $X$,
the concatenated the paths corresponding to the
$\alpha_i$ to form a piecewise Cayley graph geodesic
$\alpha_{X}$ from $x_0$ to $g\cdot x_0$.

\medskip

\noindent {\bf Claim:}  Each vertex of $\alpha_{X}$ is within
$L$ of a point of $r$ where $L$ depends only on $K$ and the quasi-isometry
constants coming from the natural quasi-isometry between $G$ and $X$ guaranteed by
Theorem ~\ref{milnor}.

Indeed, suppose $v\cdot x_0$ is a vertex on $\alpha_{X}$. Then
$v$ is either one of the $g_i$'s or $v$ lies on $\alpha_i$ for some $i$.
In the first case, $v\cdot x_0$ is within $K$ of $r$ by construction, in
the second case, $v\cdot x_0$ is within $\min\{d(v\cdot x_0,g_i\cdot
x_0),d(v\cdot x_0,g_{i+1}\cdot x_0)\} + K$ of $r$.  It suffices to
show that the first term of this sum is bounded by an $L$ as in the
statement of the claim.

This follows from the fact that the quasi-isometry between $G$ and
$X$ is proper - i.e. the ball of radius $K$ in $X$ is mapped to a
ball in $G$ whose radius is linearly distorted by the
quasi-isometry.  The vertex $v$ is on a path of the form
$g_{i+1}^{-1}g_i$ in $G$ and all of these are contained in the
pre-image of the ball of radius $2K+1$ in $X$ under the
quasi-isometry by assumption.  This gives a bound in $G$ on the
distance $v$ can be from the closest $g_i$, but now pushing back
to $X$ under the quasi-isometry gives the desired result.

Thus we now have a Cayley graph path $\alpha_{X}$ from $x_0$ to
$g\cdot x_0$ which lies in the $L$ neighborhood of $r$.  We
straighten this path to a geodesic which $\delta$-tracks $r$.
The proof is by induction on $|S|=n$.

Suppose there is an $s\in S$ which occurs in $\alpha$ and the
corresponding $s$-wall $W_s$ is crossed by $\alpha$ more than
once. If there is no such $s$, then $\alpha$ is geodesic
already.  For a fixed generator $s$, there can be more than one
$s$-wall which has multiple intersections with $\alpha$, but
these can be handled independently of one another as we point out
in Remark 1 below.  Start by finding the first occurance of $s$
along $\alpha$ so that the corresponding $s$-wall $W_s$ crosses
$\alpha$ more than once.  We consider the two cases of whether the
number of times crossed is even or odd.

\smallskip

\noindent {\it Even case:}  In this case, $\alpha$ can be written as
$asbsc$ with the
following properties:
\begin{enumerate}
\item The $s$ edge following $a$ is the first time $\alpha$ crosses
the
$s$-wall $W_s={\text Fix}(asa^{-1})$.
\item The $s$ edge following $b$ is the last time $\alpha$ crosses
$W_s$.
\item The $b$ path is all on one side of $W_s$.
\item The paths $a$ and $c$ do not cross $W_s$ at all.
\end{enumerate}

Let $M=\max\{d(x_0,s\cdot x_0)\ :\ s\in S\}$ and consider the
$\cat(0)$ segment $[a\cdot x_0, asbs\cdot x_0]$ . The endpoints
of this segment lie in $N_M({\text Fix}(asa^{-1}))$ which is a
convex subset of $X$, therefore the entire segment lies in
there.   By Theorem ~\ref{kim},
$\langle lk(s)\rangle$ acts geometrically on this $\cat(0)$ subset of $X$. Now
the induction hypothesis applies since $\langle lk(s)\rangle$ is a
right-angled Coxeter group with no more than $n-1$ generators
($s$ is not in this subgroup), that is isometrically embedded as
a convex subset of $\Lambda$ in the word metric. Thus there is a
$\delta_{n-1}$ and a Cayley graph geodesic $\beta$ in $\langle lk(s)\rangle$
so that $\beta$ $\delta_{n-1}$-tracks $[a\cdot x_0,asbs\cdot
x_0]$ and $\beta$ is also a geodesic in $G$. Thus we can replace
the path $sbs$ with the geodesic $\beta$.

\smallskip

\noindent {\it Odd case:}  In this case, $\alpha$ can be written
in the same form as above but properties 2 and 4 do not hold.  In
particular, $1_G$ and $g$ are on opposite sides of ${\text
Fix}(asa^{-1})$ and so there is exactly one more place that
$\alpha$ crosses this wall and it occurs in the $c$ path.  More
specifically, the $s$ following the $b$ in $\alpha$ is the second
to last crossing of this wall.  We do the same procedure as above
to replace $b$ with a geodesic $\beta$ which has no more
occurances of $s$ in it.  The new path $a\beta c$ obtained crosses
${\text Fix}(asa^{-1})$ exactly once.

\medskip

\noindent{\it Remark 1.}  Note that the subpath $a$ could contain
the letter $s$, but the corresponding $s$-wall is different from
the $s$-wall used above.  By assumption, any $s$-wall that $a$
may cross can only be crossed once. Also, the subpath $c$ could
contain $s$ and could cross the corresponding $s$-wall more than
once. If so, we can do the procedure described above to this
subpath independent of the procedure done for $b$ - i.e. any $s$
deletions that occur in $c$ will occur within $c$ and not with any
previous $s$ occurances by choice of $b$.  Thus with one use of
induction, we can replace pieces of $\alpha$ to obtain a path
which has the property that any $s$-wall it crosses, it crosses
exactly once.
\medskip

\noindent{\it Remark 2.}  Notice that every point of $[a\cdot
x_0,asbs\cdot x_0]$ is
within $L$ of $r$ by the $\cat(0)$ inequality.

Do this replacement procedure for each $s\in S$ to obtain a path
$\alpha_r$
which has the property that any wall it crosses, it crosses exactly
once
- i.e. $\alpha_r$ is a Cayley graph geodesic.  To see $\alpha_r$
$\delta$-tracks $r$, we use the remarks above to conclude that
$\delta_n =n\delta_{n-1}+L$ suffices.  Indeed, we will have to
use induction at most $n$ times, once for each $s\in S$.

Now suppose $y$ is any point in $X$.  It is clear that we can find a
Cayley
graph geodesic $\alpha$ which $(\delta=\delta_n+K)$-tracks $[x_0,y]$
since
we can choose an orbit point $K$ close to $y$ and do the above
procedure.

If $r:[0,\infty)\to X$ is a geodesic ray based at $x_0$, then build a
Cayley graph geodesic ray $\alpha_r$ which $\delta$-tracks $r$ as
follows:

For each $t\in [0,\infty)$, we can find $\alpha_t$ which
$\delta$-tracks
$[x_0,r(t)]$ by the previous step.  To build $\alpha$, use the local
finiteness of $\Lambda$ - indeed, infinitely many of the $\alpha_t$
share
the same first edge, infinitely many of these share the same second
edge,
etc.  This clearly builds a geodesic ray with the necessary
properties. 

It remains to show that special subgroups of $(G,S)$ are quasi-convex in $X$. Let $a$ be an element of $\langle A\rangle$ for $A$ a subset of $S$. Then our Cayley geodesic that $\delta$ tracks $[x_0, ax_0]$ is a Cayley geodesic from $1_G$ to $a\in A$. This geodesic uses only letters of $A$ by lemma \ref{convex}. Each point of $[x_0, ax_0]$ is within $\delta$ of a vertex of our Cayley path in $X$ and each vertex of our Cayley path is in $\langle A\rangle x_0$.
\end{proof}

Continuing with the quasi-isometry of the Cayley graph into the
$\cat(0)$ space for our right-angled Coxeter group $G$ as
before we give two improved versions of Proposition ~\ref{tracking}. The 
last proposition is the version we will need in the proof of the main theorem,
however we need the following intermediate step to obtain that result. 

\begin{proposition}\label{better}
Suppose $g$ and $h$ are at distance $n$ apart in the Cayley graph
and $\alpha$ is a Cayley graph geodesic from $1_G$ to $g$. Then
there exist Cayley graph geodesics $\alpha'$ and $\beta'$ from
$1_G$ to $g$ and $h$, respectively, such that each vertex of
$\alpha'$ and $\beta'$ is at most $n$ from a vertex of $\alpha$,
and $\alpha'$ and $\beta'$ have the same initial segment from
$1_G$ to $k$ for an element $k$ with Cayley graph length
$\ell(k)\geq \ell(g)-n$.
\end{proposition}

\begin{proof} We proceed by induction on $n$.  If $n=0$
then $g=h$ so $\alpha'=\beta'=\alpha$ works. Suppose $n>0$ and
take a Cayley graph geodesic $(s_1,s_2,\ldots,s_n)$ from $g$ to
$h$. Either $\ell(gs_1)$ is one more than $\ell(g)$, or else
$\ell(gs_1)$ is one less than $\ell(g)$ by lemma ~\ref{backward}.  

In the first case, denote by $\alpha_1$ the geodesic $(\alpha,s_1)$ from $1$ to $gs_1$.
Since $h$ is at distance $n-1$ from $gs_1$, by induction hypothesis there exist
$\alpha'_1$ and $\beta'_1$ within $n-1$ of $\alpha_1$ and sharing
an initial segment of length at least $\ell(gs_1)-(n-1)$.  Since
$(\alpha'_1,s_1)$ is not geodesic, by the Deletion Condition
$gs_1=us_1v$, with $u$ and $v$ represented by segments of
$\alpha'_1$ with $s_1vs_1=v$.  In fact, by lemma ~\ref{mike}, $s_1$ commutes
with each letter of $v$.  Replace the segment of
$(\alpha'_1,s_1)$ corresponding to $s_1vs_1$ with the segment of
$\alpha'_1$ corresponding to $v$ to give a geodesic $\alpha'$
from $1_G$ to $g$ within one of $(\alpha'_1,s_1)$.  Suppose $k_1$
is represented by the common initial subpath of $\alpha'_1$ and
$\beta'_1$.  If $k_1$ is given by an initial segment in $u$ then
we simply take $\beta'=\beta$.  Otherwise write $k_1=us_1w$ for
$w$ given by an initial segment of $v$.  Then $s_1ws_1=w$ and we
take $\beta'$ to be the path $\beta'_1$ with the segment
corresponding to $s_1w$ replaced by $ws_1$.  Since $uw$
corresponds to an initial segment of $\alpha'$ at most one
shorter than $k_1$, we have $\alpha'$ and $\beta'$ share a common
initial segment of length at least $\ell(g)-n$ and since
$\alpha'$ and $\beta'$ are each at most $1$ from $\alpha'_1$ they
are also at most one $n$ from $\alpha$.

In the second case where $gs_1$ is shorter than $g$, write
$g=us_1v$ for $s_1vs_1=v$ and $s_1$ commuting with each letter in
$v$ as above.  Take $\alpha_1$ to be the geodesic obtained by
replacing the segment of $(\alpha,s_1)$ corresponding to
$s_1vs_1$ by $v$, a geodesic to $gs_1$ with $gs_1$ a distance at
most $n-1$ from $h$.  By induction hypothesis we have $\alpha'_1$
and $\beta'_1$ geodesics from $1_G$ to $gs_1$ and $h$, at most
$n-1$ apart and sharing a common initial segment of length at
least $\ell(gs_1)-(n-1)=\ell(g)-n$.  Take
$\alpha'=(\alpha'_1,s_1)$ and $\beta'=\beta'_1$.  Then $\alpha'$
and $\beta'$ are at most $n$ from $\alpha$ and share the same
common initial segment as $\alpha'_1$ and $\beta'_1$ of length at
least $\ell(g)-n$.
\end{proof}

\begin{lemma}\label{best}
There exists $c$ and $d$ such that for any $r$ and $s$,
infinite geodesic rays in $X$ based at $x_0$, that are within
$\epsilon$ of each other a distance $M$ from $x_0$, there exist
Cayley graph geodesic rays $\alpha$ and $\beta$ which
$(c\epsilon+d)$ track $r$ and $s$ respectively, and which share a
common initial segment out to a distance $M-c\epsilon-d$ from
$x_0$.
\end{lemma}

\begin{proof} Appropriate $c$ and $d$ are shown below.
Suppose $r$ and $s$ are given and take $\alpha_0$ and $\beta_0$
to be Cayley graph geodesics $\delta$-tracking $r$ and $s$ for a
$\delta$ from Proposition ~\ref{tracking}. Then at a $\cat(0)$ space distance of
no more than $M+K$ from $x_0$ (where $K$ comes from the cocompactness
of the action) we have orbit points $gx_0$ and
$hx_0$ on these a Cayley graph distance $n$ apart with $n$ bounded
by a function of $\epsilon$ determined by the quasi-isometry
constants and $\delta_0$.  By Proposition ~\ref{better}, we can replace the
segments of $\alpha_0$ and $\beta_0$ out to $gx_0$ and $hx_0$ by
Cayley graph geodesic segments agreeing to a point within $n$ of
$gx_0$ in the Cayley graph, a distance in the $\cat(0)$ space
bounded again in terms of quasi-isometry constants to be at most
$(c\epsilon+d)$ from $gx_0$ and hence at least $(M-c\epsilon-d)$ from
$x_0$.  Since we are replacing an initial geodesic segment of
$\alpha_0$ by a geodesic out to $gx_0$ we again get an infinite
geodesic ray $\alpha$, actually unchanged after $gx_0$, and
similarly we get a $\beta$ from $\beta_0$. Since $\alpha$ and
$\alpha_0$ are within a distance $n$ in the Cayley graph, by a
constant determined from the quasi-isometry constants and
$\delta$, the new $\alpha$ tracks the original $r$ within some
$(c\epsilon+d)$ (we may as well take the larger of the $c$ and $d$
from here and before) and similarly for $\beta$ tracking $s$.
\end{proof}

\section{The Filtering Process}

Recall that our ultimate goal is to start with two geodesic rays in the $\cat(0)$ space $X$
whose endpoints are close in $\partial_{\infty} X$ and build a small connected set
containing these two endpoints.  From the two rays in $X$, we first obtain two 
geodesic rays $\alpha$ and $\beta$ in $\Lambda$ with a long common initial segment using lemma ~\ref{best}.   This is the starting point for the work in this section.  

Our first goal is to use the rays $\alpha$ and $\beta$ to build a collection of graphs in the plane whose edges are labeled by elements of $S$.  This collection will consist of overlapping {\it fans} (see definition ~\ref{fan}) from which we will extract a planar, one-ended graph called a {\it filter} for $\alpha$ and $\beta$.  One should think of this process of creating a filter as trying to fill in the Cayley graph between $\alpha$ and $\beta$ in a systematic way so as to avoid product subgroups yet also building a sequence of geodesics between $\alpha$ and $\beta$.  

There is an obvious map from any such filter into $\Lambda$ and this map into $\Lambda$ will be proper.  A continuous map $f:X\to Y$ is {\it proper} if for each compact set $C\subset Y$,
$f^{-1}(C)$ is compact in $X$.   The most important property of this filter is that when 
viewed in $\Lambda$, there are no long {\it factor paths} (see definition ~\ref{factor} below).  These 
are exactly the paths that behave badly under quasi-isometry so we must control the size of 
these paths if we hope to get any information to the $\cat(0)$ space $X$ from this filter.  
Moving this information into $X$ will ultimately enable us to achieve our
goal of constructing a ``small'' connected set in $\partial_{\infty} X$
containing the two points, thus showing $X$ has
locally connected boundary.  

Suppose $(G,S)$ is a right-angled Coxeter system, $G$ is a one ended group and
$G\not = \langle A\rangle \oplus \langle B\rangle$ where $A$ and
$B$ are non-empty disjoint subsets of $S$ such that $A\cup B=S$.
Recall that when $G$ is right-angled, $lk^2(C)$ where $C$ is a subset of $S$ is the ordinary link, $lk(C)$, in $\Gamma$ (see remarks following definition ~\ref{link}).  We define the notion of {\it product separator} below.  This should be compared with definition ~\ref{VFS} from earlier in the paper as both will be needed here. 

\begin{defn}[Product separator] Suppose $(G,S)$ is a Coxeter system and $\Gamma$ is the presentation graph
for $(G,S)$. A {\it product separator} for $\Gamma$ is a subset
$C$ of $S$ where $C$ is a full subgraph, $C=A\cup B$, $\langle
C\rangle=\langle A\cup B\rangle=\langle A\rangle\oplus\langle
B\rangle$, both $\langle A\rangle$ and $\langle B\rangle$ are
infinite and $C$ separates $\Gamma$.
\end{defn}

\begin{defn}[Factor path]\label{factor} Call a path $(c_{1},\ldots ,c_{m})$ in
$\Lambda$ a {\it factor path} if $\{ c_{1},\ldots ,c_{m}\}\subset
 A \cup  B$ where $A$ and $B$ are
disjoint commuting subsets of $S$ and $\langle A\rangle $ and $\langle
B\rangle $ are infinite.
\end{defn}

\noindent{\bf Remark.}\label{one-ended} Suppose $(G,S)$ is a Coxeter system with
presentation graph $\Gamma$. $G$ is one-ended if and only if $\Gamma$
contains no complete separating subgraph, the vertices of which
generate a finite subgroup.  A proof of this can be found in
\cite{MihalikTschantz}.

The following lemma allows us to easily identify when the presentation graph $\Gamma$ for $G$
has a virtual factor separator given that there are no product separators.  This lemma will be
used repeatedly in lemmas ~\ref{infinite} and ~\ref{finite} below.

\begin{lemma} \label{graph}Suppose $\Gamma$ has no product separators. Then
$\Gamma$ has a virtual factor separator iff there is a vertex $x\in
\Gamma$ such that $lk(x)$ is a subset of a visual product with
infinite factors.
\end{lemma}

\begin{proof} Suppose $\langle A\rangle \oplus \langle B\rangle$ is a
visual product in $(G,S)$ with infinite factors,
$\langle B\rangle $ has finite index in $\langle B'\rangle$ and $B'$
separates $\Gamma$ (i.e. $(B',B,A)$ is a virtual factor separator).

By the Finite Index Lemma (lemma ~\ref{finiteind}), there exists $B_{1}\subset B$
such that $\langle B'-B_{1}\rangle $ is finite and commutes with $B_{1}$. Hence
$\langle B'\rangle=\langle B_{1}\rangle \oplus \langle
B'-B_{1}\rangle$.

Note that $A\cup B'\not =S$, as $A\cup (B'-B_{1})$ commutes with
$B_{1}$ and $G$ does not decompose as a non-trivial visual direct product.
Let $y$ be a vertex in $\Gamma -(A\cup B')$ and $x$ be a vertex of
$\Gamma$ separated from $y$ by $B'$. We have $x\in A$, for
otherwise, $\langle A\cup B'\rangle =\langle A\cup (B'-B_{1})
\rangle\oplus \langle  B_{1}\rangle$ is a product
separator for $\Gamma$. If $c\in lk(x)$, and $c\not\in B'$, then $x$
and $c$ are in the same component of $\Gamma-B'$. Hence $B'$
separates $c$ and $y$ and as above $c\in A$. This means that
$lk(x)\subset A\cup B'\subset \langle A\cup (B'-B)\rangle \oplus
\langle B_{1}\rangle$.

To see the converse, suppose $x$ is a vertex of $\Gamma$ and
$lk(x)\subset (A\cup B)\subset S$ where $A$ commutes with $B$ and
both $\langle A\rangle$ and $\langle B\rangle$ are infinite. The
group $\langle lk(x)\rangle $ is infinite as $G$ has one end. If
both $\langle lk(x)\cap A\rangle$ and $\langle lk(x)\cap B\rangle$
were infinite, then $\Gamma$ has a product separator. If one, say $\langle
lk(x)\cap A\rangle$, is finite, then $\langle lk(x)\cap B\rangle$
commutes with $A$  and so $(lk(x),lk(x)\cap B,A)$
is a virtual factor separator. 
\end{proof}

\noindent{\bf Remark.} For the rest of this section, in addition
to the basic assumptions at the beginning of the section, assume
$\Gamma$ has no product separator and no virtual factor separator.

\medskip

\noindent{\it Notation.} The vertices of $\Gamma$ are the
elements of $S$ and each edge of $\Lambda$ is labeled by an
element of $S$. If $e$ is an edge in $\Lambda$, then we let $\bar
e\in S\subset G$ denote the label of $e$. Denote the vertex sets
of $\Gamma$ and $\Lambda$ by $\Gamma ^{0}$ and $\Lambda ^{0}$,
respectively.  Also, recall that $B(g)$ is the set of elements of
$S$ that make $gs$ shorter than $g$ and that this generates a
finite subgroup (lemma ~\ref{backward}).

To build a fan, we start with edge paths $(e_{1},\ldots ,e_{n},a)$ and $(e_{1},\ldots ,
e_{n},b)$ beginning at a base point $\ast$, that are geodesic in $\Lambda$, and let $g\equiv \bar
e_{1}\cdots \bar e_{n}\in G$.  Suppose $(e_{i},\ldots ,e_{n})$ is  the longest
(terminal) factor path in $(e_{1},\ldots ,e_{n})$. This may in fact be a
trivial path, as in the case $G$ is word hyperbolic.

If $\langle \bar e_{i},\ldots ,
\bar e_{n}\rangle $ {\it is infinite} and $A$ and $B$ are as in the
definition
of factor path, then $\langle \lbrace \bar e_{i},\ldots ,\bar
e_{n}\rbrace \cap
A\rangle$ or $\langle \lbrace \bar e_{i},\ldots ,\bar e_{n}\rbrace
\cap
B\rangle$ is infinite. Without loss, assume $C\equiv \lbrace \bar
e_{i},\ldots
,\bar e_{n} \rbrace \cap A$ and $\langle C\rangle $ is infinite. Now
$\lbrace
\bar e_{i},\ldots , \bar e_{n}\rbrace \subset C\cup lk(C)$. Also
$B\subset
lk(C)$ so $\langle C\rangle$ and $\langle lk (C)\rangle $ are
infinite.

\begin{lemma}\label{sep} If $\langle \bar e_{i},\ldots ,
\bar e_{n}\rangle $ is infinite, then
$\lbrace \bar e_{i},\ldots ,\bar e_{n}\rbrace\cup B(g)\subset C\cup
lk(C)$. Hence $\langle C\rangle \oplus  \langle lk(C)\rangle$ is a visual
direct product with both factors infinite and so $C\cup lk(C)$ does not
separate $\Gamma$.
\end{lemma}

\begin{proof} Recall, $s\in S$ is in $B(g)$ if $gs$ is shorter than
$g$.  If
$s\in B(g)$, then let $j$ be the largest integer such that for the
path $(e_{1},\ldots
,e_{n},s)$, $s$ deletes with $e_{j}$. If
$j\geq i$, then $s\in \lbrace \bar e_{i},\ldots ,\bar e_{n}\rbrace
\subset
C\cup lk(C)$. If $j<i$, then by lemma ~\ref{mike}, $s\in lk(C)$. In any case
$B(g)\subset
C\cup lk(C)$.
\end{proof}

\noindent{\bf Important:} The next two lemmas will allow us to choose a path in $\Gamma$ from $a$ to $b$ that avoids certain subsets of $\Gamma$.  We use this path to build a fan between these two edge paths.  We split this into two lemmas depending on whether $\langle \bar e_{i},\ldots ,\bar e_{n}\rangle$ is finite or infinite.  Note that we allow the case $a=b$.

\begin{lemma}\label{infinite} If $\langle \bar e_{i},\ldots ,
\bar e_{n}\rangle $ is infinite, then there is an edge path
$\tau$ in $\Gamma$, from $a$ to $b$ of length at least 2, such
that other than possibly a and b, no vertex of this path is in
$C\cup lk(C)$. In particular, no vertex of $\tau$ is in $\{ \bar e_i,\ldots ,\bar e_n\}$.
\end{lemma}
\begin{proof}
{\bf Case 1.} Suppose $a=b$.

By Lemma ~\ref{graph}, there is a vertex $v\in lk(a)$ such that $v\not \in
C\cup lk(C)$.
If $e$ is the edge of $\Gamma$ from $a$ to $v$ then our path is
$(e,e^{-1})$.

{\bf Case 2.} Suppose $a\in C\cup lk(C)$ and $b\not\in C\cup lk(C)$
(or when the roles of $a$ and $b$ are reversed).

By Lemma ~\ref{graph}, there is an edge from $a$ to $v\not\in C\cup lk(C)$. If
$v\not =b$, then (as $C\cup lk(C)$
cannot separate $\Gamma$) take as the desired path, the edge from $a$
to $v$ followed by any edge path from $v$ to $b$ in $\Gamma -(C\cup
lk(C))$. If $v=b$, then by Lemma ~\ref{graph} there is a vertex $w\in lk(v)$
such that
$w\not \in C\cup lk(C)$. If $e$ is the edge of $\Gamma$ from $v$ to
$w$, then take the desired path to be the edge from $a$ to $v$
followed by $e$, followed by $e^{-1}$.

{\bf Case 3.} Suppose $a,b\in C\cup lk(C)$.

By Lemma ~\ref{graph}, there is an edge $e$ from $a$ to $v\not\in C\cup lk(C)$
and an edge $d$ from $b$ to $w\not\in C\cup lk(C)$. Take as the
desired path, the edge $e$ followed by a path in $\Gamma -(C\cup
lk(C))$ from $v$ to $w$ followed by $d^{-1}$.

{\bf Case 4.} Suppose $a,b\not \in C\cup lk(C)$ and $a\not =b$.

Choose as desired path, any path of length at least 2 in $\Gamma
-(C\cup lk(C))$ from $a$ to $b$. (Again, if there is an edge between $a$ and
$b$ one can adjust as in Case 2.)
\end{proof}

\noindent{\bf Remark.} It seems somewhat artificial to consider an
edge followed by its inverse in our above paths, but paths of length at
least 2 induce beneficial combinatorics in our constructions.
These combinatorics allow for a simplified proof of lemma ~\ref{fp} at the end of 
the section.

\begin{lemma}\label{finite} If $\langle \bar e_{i},\ldots , \bar e_{n}\rangle $ is
finite, there is a path $\tau$, in $\Gamma$ from $a$ to $b$, of
length at least 2 which avoids $B(g)$.
\end{lemma}

\begin{proof} 
Recall that $B(g)$ cannot separate $\Gamma$ as $G$ is
1-ended (see remark preceding lemma ~\ref{graph}). Now proceed as in
lemma ~\ref{infinite}, case 4. 
\end{proof}

Given any edge paths $(e_{1},\ldots , e_{n},a)$ and $(e_1,\ldots ,e_n,b)$, we now have the 
path $\tau$ in $\Gamma$ from $a$ to $b$ which avoids the appropriate subsets of $\Gamma$ depending on whether $\langle \bar e_{i},\ldots ,\bar e_{n}\rangle$ is finite or infinite.  Let the ordered vertices of $\tau $ be $a=t_{0}, t_{1},\ldots,t_{m}=b$ ($m\geq 2$). Note that in $G$, $[t_{i},t_{i+1}]=1$
for all $i$.

We now describe the construction of a {\it fan} using the path $\tau$. These fans will be the buliding blocks in the construction of a {\it filter} for two geodesic rays $\alpha$ and $\beta$ in $\Lambda$.  

\begin{defn}[Fan]\label{fan} The $\tau $-{\it fan} for $(e_{1},\ldots e_{n},a)$ and $(e_{1},\ldots
e_{n},b)$ is a planar diagram consisting of edges labeled
$e_{1},\ldots , e_{n},a, b$ plus loops labeled by the vertices of
the path $\tau$ defined in Lemma ~\ref{infinite} or ~\ref{finite} (see Figure 1 below). 
The {\it fan-loops} are four sided loops
with edge labels $(t_{i}, t_{i+1},t_{i},t_{i+1})$ beginning at a
vertex labeled $g\equiv \bar e_{1}\cdots \bar e_{n}$. These loops
correspond to the commutation relations $[t_{i},t_{i+1}]$ among
consecutive vertices of  $\tau$.  We will simply refer to this as a {\it fan} for $(e_{1},\ldots ,e_{n},a)$ and $(e_{1},\ldots ,e_{n},b)$ if the path $\tau$ is unimportant. 
\end{defn}

\medskip
\begin{center}
\setlength{\unitlength}{3947sp}%
\begingroup\makeatletter\ifx\SetFigFont\undefined%
\gdef\SetFigFont#1#2#3#4#5{%
  \reset@font\fontsize{#1}{#2pt}%
  \fontfamily{#3}\fontseries{#4}\fontshape{#5}%
  \selectfont}%
\fi\endgroup%
\begin{picture}(2250,2712)(4426,-8836)
\thinlines
\put(5476,-8761){\line( 0, 1){1200}}
\put(5476,-7561){\line(-5, 4){750}}
\put(4726,-6961){\line( 0, 1){600}}
\put(4726,-6361){\line( 6,-5){450}}
\put(5176,-6736){\line( 1,-3){265.500}}
\put(5476,-7561){\line( 1, 6){150}}
\put(5626,-6661){\line(-3, 5){301.029}}
\put(5326,-6152){\line(-1,-4){148}}
\put(4951,-7386){\makebox(0,0)[lb]{\smash{\SetFigFont{12}{14.4}{\rmdefault}{\mddefault}{\updefault}$a$}}}
\put(4476,-6661){\makebox(0,0)[lb]{\smash{\SetFigFont{12}{14.4}{\rmdefault}{\mddefault}{\updefault}$t_1$}}}
\put(4951,-6496){\makebox(0,0)[lb]{\smash{\SetFigFont{12}{14.4}{\rmdefault}{\mddefault}{\updefault}$a$}}}
\put(5100,-7036){\makebox(0,0)[lb]{\smash{\SetFigFont{12}{14.4}{\rmdefault}{\mddefault}{\updefault}$t_1$}}}
\put(5400,-6926){\makebox(0,0)[lb]{\smash{\SetFigFont{12}{14.4}{\rmdefault}{\mddefault}{\updefault}$t_2$}}}
\put(5701,-6886){\makebox(0,0)[lb]{\smash{\SetFigFont{12}{14.4}{\rmdefault}{\mddefault}{\updefault}$\cdots$}}}
\put(5300,-6606){\makebox(0,0)[lb]{\smash{\SetFigFont{12}{14.4}{\rmdefault}{\mddefault}{\updefault}$t_2$}}}
\put(5486,-6336){\makebox(0,0)[lb]{\smash{\SetFigFont{12}{14.4}{\rmdefault}{\mddefault}{\updefault}$t_1$}}}
\put(5476,-7561){\line( 4, 5){600}}
\put(6076,-6811){\line( 4, 1){519.412}}
\put(6601,-6671){\line(-1,-6){ 74}}
\put(6526,-7111){\line(-5,-2){1020.345}}
\put(6151,-7486){\makebox(0,0)[lb]{\smash{\SetFigFont{12}{14.4}{\rmdefault}{\mddefault}{\updefault}$b$}}}
\put(5626,-7711){\makebox(0,0)[lb]{\smash{\SetFigFont{12}{14.4}{\rmdefault}{\mddefault}{\updefault}$e_n$}}}
\put(6226,-6661){\makebox(0,0)[lb]{\smash{\SetFigFont{12}{14.4}{\rmdefault}{\mddefault}{\updefault}$b$}}}
\put(6676,-6886){\makebox(0,0)[lb]{\smash{\SetFigFont{12}{14.4}{\rmdefault}{\mddefault}{\updefault}$t_{m-1}$}}}
\put(5476,-8561){\makebox(0,0)[lb]{\smash{\SetFigFont{12}{14.4}{\rmdefault}{\mddefault}{\updefault}-}}}
\put(5476,-7861){\makebox(0,0)[lb]{\smash{\SetFigFont{12}{14.4}{\rmdefault}{\mddefault}{\updefault}-}}}
\put(5476,-8800){\makebox(0,0)[lb]{\smash{\SetFigFont{12}{14.4}{\rmdefault}{\mddefault}{\updefault}-}}}
\put(5626,-8686){\makebox(0,0)[lb]{\smash{\SetFigFont{12}{14.4}{\rmdefault}{\mddefault}{\updefault}$e_1$}}}
\put(5976,-7086){\makebox(0,0)[lb]{\smash{\SetFigFont{12}{14.4}{\rmdefault}{\mddefault}{\updefault}$t_{m-1}$}}}
\end{picture}
\end{center}

\begin{center} {\bf Figure 1} \end{center}
\smallskip

Recall that $t_{i}\not \in B(g)$ for all $i$. Hence,
by lemma ~\ref{convex}, the edge paths $(e_{1},\ldots ,e_{n},
t_{i},t_{i+1})$ in $\Lambda$ are geodesic.
The edges labeled $a$ and $b$ at $g$ in the
$\tau$-fan are respectively called the {\it left} and {\it right fan
edges at} $g$.  The edges labeled $t_{1},\ldots ,t_{m-1}$ at $g$ are called
{\it interior fan} edges.

Next, suppose $\alpha =(e_{1},\ldots ,e_{n},a_{1},a_{2},\ldots ) $ and
$\beta=(e_{1},\ldots e_{n},b_{1},b_{2},\ldots )$ are geodesics in
$\Lambda$.

\smallskip

\noindent{\bf The Filter Construction:}  Construct a 1-ended planar graph (see Figure 2 below) with edge labels in
$S$ as follows: First construct a fan for $(e_{1},\ldots
e_{n},a_{1})$ and
$(e_{1},\ldots ,e_{n},b_{1})$.  
Next overlap this fan with a fan for
$(e_{1},\ldots e_{n},a_{1},a_{2})\textrm{ and }(e_{1},\ldots ,
e_{n},a_{1},t_{1})$ a fan for
$(e_{1},\ldots ,e_{n},t_{1},a_{1})\textrm{ and }$

\noindent$(e_{1},\ldots ,
e_{n},t_{1},t_{2}),\ldots$, a fan for $(e_{1},\ldots
,e_{n},t_{i},t_{i-1})\textrm{ and }(e_{1},\ldots
,e_{n},t_{i},t_{i+1}),\ldots$
and finally a fan for
$(e_{1},\ldots ,e_{n},b_{1},t_{m-1})\textrm{ and }
(e_{1},\ldots ,e_{n},b_{1},b_{2})$.

\begin{center}
\setlength{\unitlength}{3947sp}%
\begingroup\makeatletter\ifx\SetFigFont\undefined%
\gdef\SetFigFont#1#2#3#4#5{%
  \reset@font\fontsize{#1}{#2pt}%
  \fontfamily{#3}\fontseries{#4}\fontshape{#5}%
  \selectfont}%
\fi\endgroup%
\begin{picture}(6174,3387)(3364,-8836)
\thinlines
\put(6226,-7561){\line( 5, 3){3253.677}}

\put(5251,-6961){\line( 1, 6){150}}
\put(5401,-6061){\line( 5,-4){667.683}}
\put(6076,-6596){\line( 1,-6){160.162}}

\put(6226,-7561){\line( 3, 5){692.059}}
\put(6926,-6406){\line( 3, 1){652.500}}
\put(7576,-6195){\line(-1,-5){125.577}}

\put(4476,-6506){\line(-1, 6){ 94.162}}
\put(4381,-5941){\line( 5,-4){367}}
\put(4726,-6211){\line( 2,-3){507.692}}

\put(5241,-6961){\line(-1, 5){210.577}}
\put(5026,-5901){\line( 1, 6){ 75}}
\put(5101,-5441){\line( 1,-2){308}}

\put(7576,-6198){\line( 5, 4){375}}
\put(7961,-5891){\line(-1,-5){ 79}}
\put(7876,-6286){\line(-5,-6){425}}

\put(7426,-6846){\line( 6, 5){774.590}}
\put(8196,-6201){\line( 5, 2){430}}
\put(8616,-6036){\line(-3,-4){125}}

\put(6080,-6603){\line( 3, 4){500}}

\put(4726,-6061){\makebox(0,0)[lb]{\smash{\SetFigFont{12}{14.4}{\rmdefault}{\mddefault}{\updefault}$\cdots$}}}
\put(6301,-6661){\makebox(0,0)[lb]{\smash{\SetFigFont{12}{14.4}{\rmdefault}{\mddefault}{\updefault}$\cdots$}}}
\put(5401,-7336){\makebox(0,0)[lb]{\smash{\SetFigFont{12}{14.4}{\rmdefault}{\mddefault}{\updefault}$a_1$}}}
\put(4426,-6811){\makebox(0,0)[lb]{\smash{\SetFigFont{12}{14.4}{\rmdefault}{\mddefault}{\updefault}$a_2$}}}
\put(6826,-7456){\makebox(0,0)[lb]{\smash{\SetFigFont{12}{14.4}{\rmdefault}{\mddefault}{\updefault}$b_1$}}}
\put(6226,-8761){\line( 0, 1){1200}}
\put(6226,-7561){\line(-5, 3){2889.706}}
\put(5401,-6586){\makebox(0,0)[lb]{\smash{\SetFigFont{12}{14.4}{\rmdefault}{\mddefault}{\updefault}$t_1$}}}
\put(8026,-6136){\makebox(0,0)[lb]{\smash{\SetFigFont{12}{14.4}{\rmdefault}{\mddefault}{\updefault}$\cdots$}}}
\put(6781,-6811){\makebox(0,0)[lb]{\smash{\SetFigFont{12}{14.4}{\rmdefault}{\mddefault}{\updefault}$t_{m-1}$}}}
\put(6976,-6241){\makebox(0,0)[lb]{\smash{\SetFigFont{12}{14.4}{\rmdefault}{\mddefault}{\updefault}$b_1$}}}
\put(5886,-6916){\makebox(0,0)[lb]{\smash{\SetFigFont{12}{14.4}{\rmdefault}{\mddefault}{\updefault}$t_1$}}}

\put(6170,-6211){\makebox(0,0)[lb]{\smash{\SetFigFont{12}{14.4}{\rmdefault}{\mddefault}{\updefault}$t_2$}}}

\put(5701,-6211){\makebox(0,0)[lb]{\smash{\SetFigFont{12}{14.4}{\rmdefault}{\mddefault}{\updefault}$a_1$}}}
\put(6226,-8800){\makebox(0,0)[lb]{\smash{\SetFigFont{12}{14.4}{\rmdefault}{\mddefault}{\updefault}-}}}
\put(6226,-8536){\makebox(0,0)[lb]{\smash{\SetFigFont{12}{14.4}{\rmdefault}{\mddefault}{\updefault}-}}}
\put(6226,-7861){\makebox(0,0)[lb]{\smash{\SetFigFont{12}{14.4}{\rmdefault}{\mddefault}{\updefault}-}}}
\put(6376,-8686){\makebox(0,0)[lb]{\smash{\SetFigFont{12}{14.4}{\rmdefault}{\mddefault}{\updefault}$e_1$}}}
\put(6376,-7711){\makebox(0,0)[lb]{\smash{\SetFigFont{12}{14.4}{\rmdefault}{\mddefault}{\updefault}$e_n$}}}
\put(7876,-6886){\makebox(0,0)[lb]{\smash{\SetFigFont{12}{14.4}{\rmdefault}{\mddefault}{\updefault}$b_2$}}}
\end{picture}
\end{center}

\begin{center} {\bf Figure 2}\end{center} \smallskip

This gives our $1^{st}$ and $2^{nd}$-level fans.

At this point there is the potential for two edges of our planar
graph to share a vertex and have the same label. These edges
are {\it not} identified in our planar graph. Rather, an edge path
$\tau$, constructed as in case 1 of lemma ~\ref{infinite} will be used to extend
our graph between these two edges.

To continue, we must specify geodesics from the base point $\ast$ to each vertex
defined so far.  In each fan-loop constructed up to this point,
designate the upper left edge as a {\it non-tree} edge. The graph minus
the non-tree edges is a tree. Take as designated geodesic from $\ast$
to a defined vertex, the unique geodesic of the tree.

Continuing, (at each stage designating upper left edges of fan loops
as non-tree edges and constructing fans) we obtain the desired 1-ended planar 
graph called a {\it filter} for
$\alpha$ and $\beta$. Geodesics from $\ast$ in this
filter determine geodesics in $\Lambda$ with the same edge labels.
Hence there is a natural proper map from this filter to $\Lambda$. The
image of this filter in $\Lambda$ is a closed 1-ended subgraph.

The following are useful combinatorial facts about a filter $F$.
\begin{enumerate}
\item Each vertex of $F$ has exactly one or two
edges directly below it.  This fact would be missing if certain
identifications had been allowed.
\item If a vertex of $F$ has exactly one edge below it, then
this edge is either an interior fan edge, an $\alpha $-edge or a
$\beta $-edge.
\item If a vertex of $F$ has exactly two edges $e$ and $d$ below it, then
one (to the right) is a left fan edge, say of the fan $\tau $ and the other (to the left) is a right fan
edge of the fan $\lambda$, where $\tau\ne \lambda$ and $\bar d\ne \bar e$. Both $e$ and $d$ belong to a single fan loop.
\item For each vertex of $F$ there is a unique fan containing all edges directly above this vertex. In particular, the edges above this vertex consist of a left and right fan edge and at least one interior edge.\item The filter for $\alpha$ and $\beta$ minus all non-tree
edges is a tree containing $\alpha$ and $\beta$ and all interior edges of all fans.
\item An edge of the filter is a non-tree edge iff it is a
right fan edge not on $\beta$. Right fan edges not on $\beta$ are upper left edges of a fan-loop.
\item Let $T$ be the tree obtained from a filter $F$ by removing all
non-tree edges.  There are no dead ends in $T$ - i.e. for any vertex
$v\in T$, there is an interior edge at $v$.
\end{enumerate}

We now begin the process of showing that if we create a filter for geodesic rays $\alpha$ and $\beta$ that share a long common initial segment, then there is a bound on the length of any factor path in this filter.  Again, this will allow us to map the one-ended planar graph into $X$ and obtain a small connected set in $\partial_{\infty} X$.  The final result necessary is lemma ~\ref{mainlem} but we prove two intermediate lemmas about factor paths in filters along the way.

\begin{lemma}\label{fp} Let $\gamma\equiv (c_{1},\ldots ,c_{n})$ be an edge
path in $T$ such that each edge is above the previous edge and
at most, the initial point of $\gamma$ intersects $\alpha \cup \beta$.
Suppose $(c_{i},\ldots , c_{j})$ is a factor path
of $\gamma$ of length $> \vert S\vert $. Then,
either $\bar c_{j+1}$ or $\bar c_{j+2}$
is not in $\lbrace \bar c_{i},\ldots ,\bar c_{j}\rbrace$.
\end{lemma}

\begin{proof}
Since $j-i\geq S$, and $G$ is right angled, $\langle \bar c_i,\ldots \bar c_j\rangle$ is infinite. Neither, $c_{j+1}$ nor $c_{j+2}$ is a right fan edge, as
both are in $T$. If $c_{j+1}$ is an interior fan edge we are done by
the construction in lemma ~\ref{infinite}. Hence
we may assume that $c_{j+1}$ is a left fan edge and $\bar c_{j+1}\in
\lbrace \bar c_{i},\ldots, \bar c_{j}\rbrace$. Now, if $c_{j+2}$ is an interior
fan edge, once again we are finished. Hence we may assume $c_{j+2}$
is a left fan edge. 

Let $u$ be the end point of $c_{j+2}$ and $x$ the initial point of $c_i$. As $c_{j+2}$ is a left fan
edge, there is a right fan edge $a$ below $u$ (see (2) and (3) of the above list of combinatorial facts).

Note that, $\bar a\ne \bar c_{j+1}$ since otherwise there is no right fan edge below the end point of (the left fan edge) $c_{j+1}$.
Hence the situation must be one of the following two cases pictured in Figure 3.

{\bf Case 1.} Suppose $\bar a\not\in \lbrace \bar c_{i},\ldots ,\bar
c_{j}\rbrace$.

Consider the path $(c_{i},\ldots ,c_{j+2})$ at $z$. By (4), the edge
$c_{i}$ at $z$ is either an internal or left fan
edge, and all other edges of this path are internal fan edges. Hence
this entire path is in $T$. But then by the lemma \ref{infinite} construction,
$\bar c_{j+1}\not\in \lbrace \bar c_{i},\ldots ,
\bar c_{j}\rbrace $, contrary to our assumption.$\bullet$

\medskip

{\bf Case 2.} Suppose $\bar a\in \lbrace \bar c_{i},\ldots ,\bar
c_{j}\rbrace$.

Let $k$ be the largest integer in $\lbrace i,\ldots ,j\rbrace$ such
that $\bar a=\bar c_{k}$.  Note that $\bar a\not
=\bar c_{j+1}$ as this
was considered above . The edge path $(c_{i},\ldots ,
c_{k-1},c_{k+1},\ldots ,c_{j+2})$ at $x$ is in $T$, as (4) implies
the first edge of the subsegment $(c_{k+1},\ldots ,c_{j+2})$ (at $y$) is
either an internal or left fan edge and each subsequent edge is an internal fan edge. 
Hence $(c_{i},\ldots ,c_{k-1},c_{k+1},\ldots ,c_{j+2})$ at $x$ defines a (geodesic)
factor path. As the length of this geodesic is $\geq \vert S\vert +2$,
$\langle \bar c_{i},\ldots ,\bar c_{k-1},\bar c_{k+1},\ldots , \bar
c_{j}\rangle$ is
infinite. We have $\bar c_{j+2}\not \in \lbrace \bar c_{i},\ldots
\bar c_{k-1},\bar c_{k+1},\ldots ,\bar c_{j}\rbrace$, by our $\tau$-fan construction in lemma ~\ref{infinite}. But, we also have
$\bar c_{j+2}\not=\bar c_{k}=\bar a$. 
\end{proof}

\begin{figure}
\vbox to 3in{\vspace {-.1in} \hspace {.2in}
%\hspace{-1 in}
\includegraphics[scale=1]{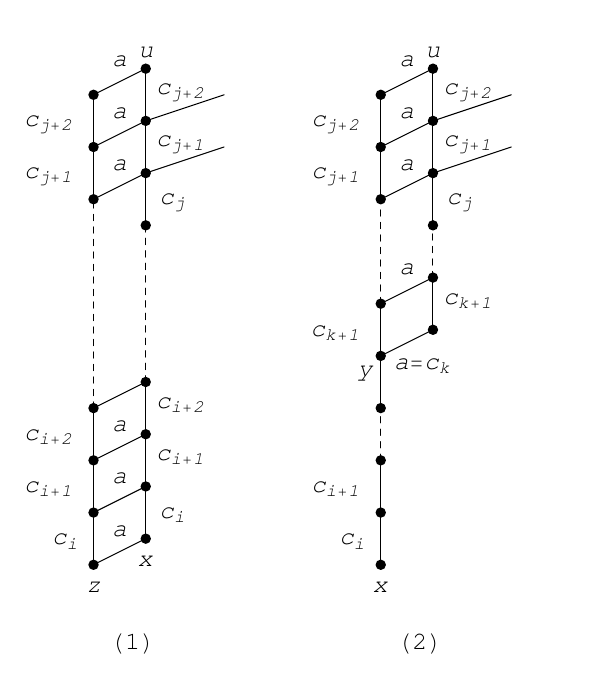}
\vss }

\vspace{1.5in}
\centerline {\bf Figure 3}
%\caption{Avoiding compact sets} 
\label{Fig231d}

\end{figure}

\begin{lemma}\label{final} Suppose that $\gamma\equiv (c_{1},\ldots c_{n})$ is a
directed factor path in $T$ such that $\gamma$ intersects $\alpha
\cup \beta$ in at most, its initial point. Then $n\leq 3\vert S\vert$.
\end{lemma}

\begin{proof}
By lemma ~\ref{fp}, at least every other edge $e_{i}$
following $e_{\vert S\vert}$ is such that $\bar e_{i}\not \in \lbrace
\bar e_{1},\ldots \bar e_{i-1}\rbrace$.
\end{proof}

The following lemma is the main lemma of this section that will be used in the 
proof of main theorem.   In order to obtain the necessary bounds on the lengths of
factor paths, we must analyze how such a path arises in $\Lambda$.  This involves studying
Van-Kampen diagrams in the group.  We refer the reader to \cite{LyndonSchupp} for the necessary background concerning Van-Kampen diagrams for groups (abbreviated V-K diagrams).

\newpage

\begin{lemma}\label{mainlem}
For any $k$, there exists a $\delta$ such that, for any $\alpha$
and $\beta$ Cayley graph geodesics between the same pair of
points, if for some $n$, $\alpha(n)$ and $\beta(n)$ are at least
$\delta$ apart, then the initial (terminal) segment of $\alpha$ that ends (begins) at $\alpha(n)$ contains a factor subpath of length $k$
(similarly for  $\beta$).
\end{lemma}

\begin{proof}Write $\alpha=(\alpha_0,\alpha_1)$ and
$\beta=(\beta_0,\beta_1)$ for $\alpha_0$ and $\beta_0$ initial
subpaths of length $n$. Take a Cayley graph geodesic $\gamma$ from
$\alpha(n)$ to $\beta(n)$ and V-K diagrams for the loops
$\alpha_0\gamma\beta_0^{-1}$ and $\gamma\beta_1\alpha_1^{-1}$ each with a
minimum number of regions, such that the combination is a V-K diagram
for the relation telling us that $\alpha$ and $\beta$ represent the same
element of the Coxeter group.  Given any edge in this V-K
diagram, we construct a chain of boxes corresponding to commuting
generators in the defining presentation.  Such a chain cannot
cross $\gamma$ twice or begin and end on $\alpha$ or begin or end
on $\beta$ since each of these is a geodesic. Such a chain cannot
be a loop else one of the diagrams wouldn't have had a minimum
number of regions.  No two chains with the same letter cross each
other.  Thus each chain begins on $\alpha$, ends on $\beta$, and
perhaps crosses $\gamma$ with as many chains crossing from
$\alpha_0$ to $\beta_1$ as from $\alpha_1$ to $\beta_0$. Note that
$\gamma$ consists of an edge from each such chain. Each letter for a
chain from $\alpha_0$ to $\beta_1$ commutes with a letter on a
chain from $\alpha_1$ to $\beta_0$ (since these chains cross at
some point) which gives a commutation relation between these letters.
Write $\gamma_0$ for the subsequence of letters of $\gamma$ that
label chains extending from $\alpha_0$ to $\beta_1$ and
$\gamma_1$ for the subsequence of letters of $\gamma$ that label
chain extending from $\alpha_1$ to $\beta_0$.  Then in fact
$\gamma$ and $(\gamma_0,\gamma_1)$ would represent the same
element since each letter of $\gamma_0$ commutes with each letter
of $\gamma_1$.

We take $\delta$ greater than twice the longest geodesic in a finite
subgroup of the Coxeter group.  Then the letters in $\gamma_0$
generate an infinite subgroup as do the letters in $\gamma_1$.  Let
$\alpha'_0$ be the subsequence of letters in $\alpha_0$ belonging to
chains crossing $\gamma$ (a word equivalent to $\gamma_0^{-1}$).  A
letter of $\alpha_0$ not in $\alpha'_0$ must commute with each earlier
letter belonging to $\alpha'_0$ since the chain extending from such a
letter extends to $\beta_0$ and crosses the chains of earlier letters
in $\alpha'_0$.  Of the $\delta/2$
letters in $\alpha'_0$, one letter (say $s_1$) must occur at least $\delta/(2|S|)$
times.  Between occurrences of $s_1$ there must be a letter
not commuting with $s_1$ and hence also a letter of $\alpha'_0$, and
among all such letters between successive $s_1$ occurrences, there must be some
letter occurring at least $(\delta/(2|S|)-1)/|S|$ times say $s_2$.  We
define a subset $S'$ of the letters $s_1,s_2,\ldots,s_i$ occurring in
$\alpha'_0$, and define a sequence of subwords of $\alpha'_0$ as
follows.

First we take $s_1$ and $s_2$ in $S'$ and take the terminal subword of
$\alpha'_0$ after the later of the first occurrences of $s_1$ and
$s_2$. Now the first occurrence of $s_3$ occurs either in the first half or
the last half of this subword of $\alpha'_0$.  If the terminal segment of
$\alpha'_0$ after the first occurrence of $s_3$ in the remaining word
is longer than the word before this $s_3$, then we keep $s_3$ in $S'$ and
take the terminal segment after this $s_3$, and otherwise we omit $s_3$
from $S'$ and restrict to the $s_3$-free initial subword of the remainder.
Repeat the process for $s_4,\ldots,s_i$ - each time taking the first
occurrence of the next letter and dividing the remaining word at this
point keeping the larger half, putting the letter in $S'$ if we take a
terminal segment, and otherwise omitting it from $S'$.

\begin{figure}
\vbox to 3in{\vspace {-1.5in} \hspace {-2in}
%\hspace{-1 in}
\includegraphics[scale=1]{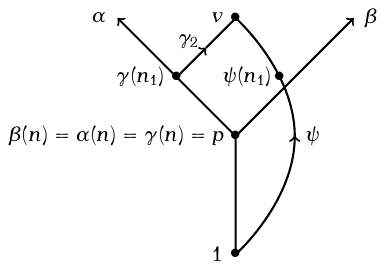}
\vss }
\vspace{-.9in}
%\caption{Avoiding compact sets} 
\label{Fig231d}

\centerline {\bf Figure 4}

\end{figure}

\medskip

In the end we are left with a subword $\alpha''_0$ of $\alpha'_0$ of
length at least $2(\delta-2|S|)/(|S|^2 2^{|S|})$ such that the following
conditions hold:  each letter in  $\alpha''_0$ belongs to $S'$ and
each letter in $S'$ occurs either at the beginning of $\alpha''_0$ or
else earlier in $\alpha'_0$ than all of $\alpha''_0$ as well.  Let $\phi$
be smallest subword of $\alpha_0$ containing the
corresponding letters in $\alpha''_0$.  Then the letters in $\phi$
belong to $S'$ or belong to the subset $lk(S')$ of letters commuting with all
elements of $S'$.  Since the letters in $\gamma_1$ are in $lk(S')$,
$lk(S')$ generates an infinite subgroup, as does $S'$ since $s_1$ and
$s_2$ are in $S'$.  Hence there is a subword in $\alpha_0$ with letters in a
product of infnite subgroups having length bounded below by a
function of $\delta$.  If $\delta\geq |S|^2 2^{|S|} k+2|S|$ then this subword
will be a factor path of length at least $k$.
\end{proof}

\section{Proof of the Main Theorem}

\begin{thm}\label{mainthm}Suppose $(G,S)$ is a right-angled Coxeter system, $G$ is 1-ended,
$\Gamma (S)$ contains no product separator and no virtual factor
separator. Also assume that $G$ does not visually split as a
non-trivial direct product. Then $G$ has locally connected boundary.
\end{thm}

\begin{proof} Suppose $G$ acts geometrically on a $\cat(0)$ space $X$.
Choose a base point $\ast \in X$. Suppose $r$ and $s:[0,\infty )\to X$
are geodesics at $\ast$ and $s\in N_{(M,\epsilon)}(r)$ (i.e.
$d_{X}(r(M),s (M))<\epsilon$).

Choose geodesics $\alpha =(e_{1},\ldots ,e_{n},a_{1},a_{2},\ldots )$ and
$\beta =(e_{1},\ldots ,e_{n}, b_{1},b_{2},\ldots )$ from $\Lambda$
that $\delta$ track $r$ and $s$ respectively as in lemma ~\ref{best}.
Recall $\alpha_X$ means view the path $\alpha$ inside $X$.  Thus we can
choose $\alpha$ and
$\beta$ so that $\alpha_X$ and $\beta_X$ share a common initial piece that is
$M'$ long where $M'$
comes from following the $M-c\epsilon -d$ from lemma ~\ref{best} under the
quasi-isometry between $\Lambda$ and $X$.  Note that $d(\alpha_X(M'),
r(M'))<\delta $ and $d(\beta_X(M'),r(M'))<\delta $.

Construct a filter $F$ for this pair of edge paths and consider the tree
$T$ constructed from $F$ as in section 5. Let $f:F\to X$ be the natural map
that factors through $\Lambda(G,S)$. To show local connectivity at $r$, it
suffices to show the limit set of $f(F)$ is a ``small" connected subset of
$\partial X$ containing both $r$ and $s$.

The map $f$ is proper and $F$ is one-ended, thus the limit set of $f(F)$
is connected.  It is clear that this limit set contains $r$ and $s$.
To see that this limit set is small, we apply lemma ~\ref{mainlem} and lemma ~\ref{final} to show that: for any vertex $v\in F$, the CAT(0) geodesic from $\ast$ to $f(v)$ passes close to $r(M')$. If $r$ and $s$ determine points close to one another in the limit set of $X$, then $M'$ is a large number. 
Here is the important point. 

\medskip

\noindent{\it Claim:}  There exists $\epsilon '$ (independent of $M'$) such that
$f(F)\subset N_{(M',\epsilon ')}(r)$ (the $M'$ mentioned above is the
same $M'$ used here).

First we work in $F$, then transfer this information over to $X$ via map
$f$.

Let $v$ be a vertex in $F$ that is far away from $\alpha(n)$.
In particular, we want $d(\ast, f(v))$ to be much
bigger than $M'$ since we are concerned with the limit set of $F$.  Also,
denote by $p$ the point at which $\alpha$ and $\beta$ start to separate.

We know the geodesic in $F$ from the identity to $v$ is obtained by
following $\alpha$ (and $\beta$) to $p$, then using the unique geodesic in
the tree $T$ (described in section 5) from $p$ to $v$.  Call this geodesic
$\gamma$.  Thus $\gamma (n)=\alpha (n)=\beta (n)=p$. Write $\gamma=(\gamma_1,\gamma_2)$ where $\gamma_1(n_1)$ is the last vertex of $\gamma_1$ and the last vertex of $\gamma$ on $\alpha\cup \beta$. Note that $n_1\geq n$.  Lemma \ref{final} provides a bound $k$, on the length of factor subpath in $\gamma_2$.
%By lemma ~\ref{final}  there exists a global bound on the length of a (directed) factorpath in $T$ that only intersects $\alpha\cup\beta$ in it's initial point,call this bound $k$.  
By lemma ~\ref{mainlem} (applied to $\gamma_2$ and $n_1$) there is a $\delta_1$ coming from this
$k$ such that every Cayley geodesics $\psi$, from the identity to $v$, is such that $\psi(n_1)$ is contained in
the $\delta_1$ neighborhood of $\gamma (n_1)$ (a vertex of $\alpha\cup \beta$ beyond $p$.

We transfer all of this information over to $X$ via the quasi-isometry to
obtain new constants.  In particular, the $\cat(0)$ geodesic, $\tau$ from
$\ast$ to $f(v)$ is within  $\delta_2$ of $f(\gamma (n_1))$.  Thus
$d(\tau(M'), f(p))<\delta_2$.

We know from the choice of $\alpha$ and $\beta$ above that
$d(\alpha_X(M'), r(M'))<\delta $ and $d(\beta_X(M'),r(M'))<\delta $.
But $\alpha_X(M')$ is simply $f(p)$ here, thus we have
$d(f(p),r(M'))<\delta $.  Putting these together with the triangle
inequality yields:  $d(r(M'),\tau (M'))< \delta _2+\delta $.  If we use
the $M'$ mentioned above along with $\epsilon '=\delta_2+\delta$, we have
the claim. 
\end{proof}

\end{document}